\documentclass[Journal]{IEEEtran}
\IEEEoverridecommandlockouts
\ifCLASSINFOpdf
\else
\fi

\usepackage{mathrsfs}
\usepackage{color}
\usepackage{amsfonts,amssymb}
\usepackage{amsfonts}
\usepackage{subfigure}
\usepackage{booktabs}
\usepackage{graphicx} 
\usepackage{float} 
\usepackage{amsmath}
\usepackage{enumerate}
\usepackage{algorithm}
\usepackage{algpseudocode}
\usepackage{bm}
\usepackage{makecell}
\usepackage{epstopdf}
\usepackage{mathtools}
\usepackage{bbm}
\usepackage{dsfont}
\usepackage{pifont}
\usepackage{amsthm}
\usepackage{cite}
\usepackage{balance}
\usepackage[T1]{fontenc}
\usepackage{tikz}

\newtheorem{remark}{Remark}

\newtheorem{theorem}{Theorem}
\newtheorem{lemma}{Lemma}

\newtheorem{assumption}{Assumption}

\newtheorem{definition}{Definition}
\newtheorem{proposition}{Proposition}

\begin{document}
\title{\LARGE \bf Near-Optimal Mixed Strategy for Zero-Sum Differential Games}
\author{Tao Xu, Wang Xi, and Jianping He
	\thanks{The authors are with the Department of Automation, Shanghai Jiao Tong University, and Key Laboratory of System Control and Information Processing, Ministry of Education of China, Shanghai 200240, China. E-mail: \{Zerken, bddwyx, jphe\}@sjtu.edu.cn. Preliminary results of this work were presented at IEEE Conference on Decision and Control (CDC), 2023 \cite{xuDifferentialGameMixed2023b}.} }

\maketitle

\begin{abstract}
    Synthesizing near-optimal mixed strategies for zero-sum differential games (ZSDGs) has been a longstanding challenge. Existing research mainly focuses on characterizing the theoretical value function, while the practical design of executable mixed strategies remains open. To address this issue, we propose a novel weak approximation framework. The core idea is to map the original mixed-strategy game into a surrogate stochastic differential game (SDG) under pure strategies. This mapping ensures that both state distributions and cost expectations closely match the original game. Based on the solution of this auxiliary SDG, the original game value can be approximated, and near-optimal mixed strategies can be synthesized. To operationalize this framework, we develop a constructive control-space discretization algorithm for general ZSDGs. By parameterizing the infinite-dimensional measure optimization into standard probability simplices and solving local linear programs, our method efficiently synthesizes executable mixed strategies. Furthermore, we rigorously prove that the global weak approximation error is strictly of order $\mathcal{O}(\bar\pi)$ with respect to the maximum commitment delay $\bar\pi$, and derive explicit analytical upper bounds for the strategy suboptimality gaps. Numerical examples are provided to illustrate and validate our theoretical results.
\end{abstract}

\section{Introduction}
Zero-sum differential games (ZSDGs) are powerful in modeling dynamic conflicts, where one player seeks to minimize a cost functional while the opponent aims to maximize it. Originating from Isaacs’ pioneering work \cite{isaacsDifferentialGamesMathematical1999} and further refined through the Hamilton–Jacobi–Isaacs (HJI) framework \cite{evansDifferentialGamesRepresentation1984c}, these games traditionally rely on pure strategies to attain saddle-point solutions. However, when pure strategies fail to yield an equilibrium, \textit{mixed strategies}, where players randomize over their control actions, become essential. Mixed strategy is proved indispensable in various fields such as economics \cite{harrisExistenceSubgameperfectEquilibrium1995}, generative adversarial networks \cite{goodfellowGenerativeAdversarialNets2014}, reinforcement learning \cite{perkinsMixedStrategyLearningContinuous2017a}, and robotics \cite{petersLearningMixedStrategies2022}.

Historically, much of the literature has focused on characterizing the value function of ZSDGs via the HJI framework. A key concept in this context is the notion of nonanticipative strategies with delay (NAD) \cite{varaiyaExistenceSaddlePoints1969,roxinAxiomaticApproachDifferential1969,elliottExistenceValueStochastic1976}, which require players to commit control actions for a predetermined delay period, referred to as \textit{commitment delay}, before responding to the opponent’s moves \cite[p.114]{engwerdaLQDynamicOptimization2005}. This delayed commitment not only mirrors realistic constraints but also helps manage analytical challenges inherent in differential games. 

\subsection{Motivations}
Despite significant progress in characterizing the value function of ZSDGs under mixed strategies via the HJI framework, efficiently computing the optimal mixed strategies remains an open and challenging problem. In current approaches \cite{cardaliaguetDifferentialGamesAsymmetric2007,cardaliaguetPureRandomStrategies2014,buckdahnValueFunctionDifferential2013,buckdahnValueMixedStrategies2014,flemingMixedStrategiesDeterministic2017a}, the commitment delays in mixed strategies are required to asymptotically approach zero, which leads to chattering behavior for mixed strategy design—that is, the necessity of switching control inputs infinitely fast \cite{levantChatteringAnalysis2010}. This not only poses severe implementation challenges in practical applications, such as robotics and pursuit-evasion, but also limits the theoretical robustness of these solutions. Addressing this gap is critical, as developing methods to compute optimal mixed strategies without relying on infinitesimal commitment delays could significantly enhance the real-world applicability of mixed strategies in dynamic conflict scenarios.

In this paper, we aim to develop a weak approximation framework to solve ZSDGs under mixed strategies without the asymptotic assumption on the commitment delay. The key idea is to design a stochastic equivalent formulation under pure strategies that approximates the original game under mixed strategies. By ensuring that the deviations between both the state distributions and the cost expectations remain small over the entire time horizon, the proposed framework enables us to approximate both the value function and the optimal mixed strategies effectively. 

Our weak approximation framework draws inspiration from research dedicated to understanding the dynamics of stochastic gradient descent algorithms \cite{liStochasticModifiedEquations2017, liStochasticModifiedEquations2019,evenContinuizedAccelerationsDeterministic2021,latzAnalysisStochasticGradient2021a}. In these studies, the authors devise stochastic differential equations (SDEs) to approximate the distributional dynamics of stochastic gradient descent (SGD) algorithms in a weak sense. We have observed a shared challenge between understanding SGD algorithms and characterizing game dynamics and accumulated costs under mixed strategy. Building upon this similarity, the idea of employing SDE to approximate SGD algorithms can be further extended to the proposed SDG weak approximation.

\subsection{Challenges}
To solve the optimal mixed strategies for a ZSDG based on an SDG weak approximation method, the first challenge lies in a reasonable definition of admissible mixed strategy spaces. Nevertheless, a universally accepted definition of mixed strategies for ZSDGs is still absent \cite{cardaliaguetPureRandomStrategies2014}. For zero-sum static games, defining mixed strategies can usually be simplified to specifying probability distributions over the admissible pure strategy space \cite[p.32]{osborneCourseGameTheory1994}. For ZSDGs, directly defining probability distributions over the admissible pure strategy spaces leads to a dead end caused by some well-known measure-theoretic obstacles \cite{aumann28MixedBehavior1964}. To guarantee the existence of saddle-point equilibrium, it requires delicate constructions to endow the measurable structure \cite{buckdahnValueFunctionDifferential2013,buckdahnValueMixedStrategies2014,cardaliaguetPureRandomStrategies2014}. To further enable efficient solutions of optimal mixed strategies, it put forward higher requests on designing admissible mixed strategies.

The second challenge comes from the development of the SDG weak approximation. Compared to the SDE weak approximations of SGD algorithms, SDG weak approximation of ZSDGs under mixed strategies is much more complicated. On the one hand, the randomness in an SGD algorithm only comes from a pre-defined uniform distribution over a finite set of gradient indexes, while the randomness caused by mixed strategies is dynamically evolving along the time horizon over the continuous set of control spaces. On the other hand, SDE weak approximation aims to approximate the state distributions, while SDG weak approximation requires approximating not only the state distributions but also the cost expectations. Furthermore, it also needs to certify that the value function of the designed SDG under pure strategies can well approximate the original game under mixed strategies. 

The third challenge arises from the algorithmic realization of near-optimal mixed strategies. While the proposed weak approximation framework lays the theoretical foundation for matching state distributions, deriving the exact optimal probability measures remains analytically intractable. Designing near-optimal mixed strategies essentially requires bridging the gap between the infinite-dimensional measure optimization and a computationally solvable mechanism, and further certifying the suboptimality gaps of the resulting strategies when mapped back to the original game.

\subsection{Contributions}
In this paper, we develop a weak approximation framework coupled with a constructive discretization algorithm to address the aforementioned challenges. 
The main contributions are summarized as follows:
\begin{itemize}
    \item \textbf{(Existence)} By characterizing the admissible mixed strategies via stepwise functions with non-vanishing commitment delays, we rigorously establish the existence of a game value for ZSDGs. This physics-aware formulation essentially bridges the theoretical measure optimization with practical execution, inherently circumventing the un-implementable chattering behaviors in control inputs.
    
    \item \textbf{(Framework)} A novel weak approximation framework is established to map the original mixed-strategy ZSDG into a surrogate pure-strategy SDG. This theoretical paradigm simultaneously tracks the state distributions and the cost expectations along the time horizon. Furthermore, it ensures that the value function of the auxiliary SDG robustly approximates that of the original game.
    
    \item \textbf{(Algorithm)} We develop a finite-dimensional control-space discretization algorithm to design near-optimal mixed strategy for the general game. By parameterizing the infinite-dimensional measures into probability simplices and solving local linear programs, this method efficiently synthesizes executable near-optimal mixed strategies. Furthermore, we rigorously certify a global weak approximation error bounded by $\mathcal{O}(\bar\pi)$ and derive explicit analytical suboptimality gaps.
\end{itemize}

The rest of this paper is organized as follows. Sec. II reviews the related works. Sec. III introduces some preliminaries and formulates the problem of interest. In Sec. IV, the admissible mixed strategy spaces are proposed, under which the existence of the game value is proved. Sec. V expounds on our theoretical weak approximation framework in a five-step manner. Sec. VI introduces a control-space discretization algorithm for general ZSDGs. Simulation results are presented in Sec. VII. Finally, we conclude in Sec. VIII.

\section{Related Works}
In this section, we provide a brief review of the extensive research on the definitions, solutions and applications of mixed strategies for ZSDGs.

The most direct idea, following the mixed strategy of static games \cite{nashjrNoncooperativeGames1996}, is to define probability distributions over the admissible pure strategy space. However, it suffers from the measure-theoretic problem of choosing a proper subset as delineated in \cite[p.231]{basarDynamicNoncooperativeGame1998}. Rather than endowing measurable structure on the pure strategy space, it is suggested to define mixed strategies as random variables in the seminal work by R. J. Aumann \cite{aumann28MixedBehavior1964}. This remarkable idea makes the definition of mixed strategy much more flexible and thus has inspired a lot of subsequent research \cite{cardaliaguetDifferentialGamesAsymmetric2007, cardaliaguetPureRandomStrategies2014,flemingMixedStrategiesDeterministic2017a,buckdahnValueFunctionDifferential2013,buckdahnValueMixedStrategies2014}. Another line of work \cite{sirbuStochasticPerronsMethod2014,kunischOptimalControlUndamped2016,coxControlledMeasurevaluedMartingales2024} compose their definitions following the concept of relaxed controls \cite{berkovitzRelaxedControls2012}, where the admissible controls are measure-valued functions. Being more interested in real-valued control functions in this paper, we follow the ideas from Aumann. In \cite{cardaliaguetDifferentialGamesAsymmetric2007,cardaliaguetPureRandomStrategies2014,flemingMixedStrategiesDeterministic2017a}, measurable structures are endowed during the definitions of admissible strategies, and a mixed strategy is a measurable map from the sample and information set to an admissible control. In \cite{buckdahnValueFunctionDifferential2013,buckdahnValueMixedStrategies2014}, the authors choose to endow a measurable structure inside the definition of admissible control, and then define admissible mixed strategies based on the concept of stochastic interval. To guarantee the game value exists, these works require commitment delays to asymptotically approach zero.

Even for the differential game under pure strategies, the game value generically has no explicit expression or simply does not exist \cite{elliottExistenceValueStochastic1976}. Traditional methods to solve differential games under pure strategies include the Pontryagin maximum principle \cite{wangPontryaginsMaximumPrinciple2010} and viscosity solutions of HJI equations \cite{flemingControlledMarkovProcesses2006}. However, there is relatively little research on solving optimal mixed strategies for ZSDGs. In \cite{kumarOptimalMixedStrategies1980}, the optimal mixed strategies to a class of generalized bomber-and-battleship game are solved. In \cite{flemingMixedStrategiesDeterministic2017a}, the existence of approximately optimal mixed strategies for zero-sum differential games is proved. 
Another line of research is the learning-based method that parametrizes mixed strategies by neural networks with input noises, including \cite{dou2019finding} and \cite{martinFindingMixedstrategyEquilibria2023}.

The idea of randomized strategies has been widely applied to differential games such as pursuit-evasion problems \cite{hespanhaProbabilisticPursuitevasionGames2000,vidalProbabilisticPursuitevasionGames2002a, islerRandomizedPursuitevasionPolygonal2005, islerRandomizedPursuitEvasionLocal2006}. Authors of \cite{hespanhaProbabilisticPursuitevasionGames2000} consider evaders that actively avoid the pursuer, then describe the probabilistic pursuit-evasion game as a partial information Markov game and introduce the discrete stochastic strategies to solve it. \cite{vidalProbabilisticPursuitevasionGames2002a} continues to study the intelligent randomized evader strategies by further considering vision-based detection. Authors of \cite{islerRandomizedPursuitevasionPolygonal2005} investigate the randomized pursuit strategies and provide a drastic increase in the power of the pursuer even if the evader is faster than the pursuer and knows the position of the pursuers all the time. Authors of \cite{islerRandomizedPursuitEvasionLocal2006} further study the case where only local visibility is available for the pursuer. These works show that mixed strategies can improve performances for both pursuer and evader in differential games. However, equilibrium is not guaranteed for these proposed strategies because of the difficulty in optimizing over mixed strategies.

\section{Preliminaries and Problem Formulation}
\subsection{Preliminaries}
\textbf{Notations.} Let $\pi = \{t_{0} < t_{1} < \ldots < t_{N} = T\}$ be a partition of a time interval $[t_0, T]$, and $\Pi_{[t_0,T]}$ be the set of all partitions on $[t_0,T]$. We define the lower bound, upper bound, and order of $\pi$ as follows,
\begin{equation*}
    \underline{\pi} \triangleq \min_{0 \leq k \leq N-1}(t_{k+1}-t_{k}),\; \bar{\pi} \triangleq \max_{0 \leq k \leq N-1}(t_{k+1}-t_{k}),\;|\pi|\triangleq N.
\end{equation*}
A partition $\pi$ is called equispaced if $\underline{\pi} = \bar{\pi}$. 
The $k$-th time interval, $[t_{k-1}, t_k)$ is denoted as $\Delta_k$, with length $|\Delta_k|\triangleq t_k-t_{k-1}$. The indicator function $\mathbb{I}_{\Delta_k}(t)$ equals $1$ if $t\in \Delta_k$, otherwise equals $0$. 
$\zeta_\pi:[t_0,T]\to\{1,\ldots,|\pi|\}$ is the index function of $\pi$, and $\delta_\pi:[t_0,T]\to[\underline{\pi},\bar{\pi}]$ is the interval length function of $\pi$, i.e., 
\[\zeta_{\pi}(t) \triangleq \sum_{k=1}^{|\pi|} k\,\mathbb{I}_{\Delta_k}(t), \; \delta_{\pi}(t)\triangleq\sum_{k=1}^{|\pi|}|\Delta_k|\,\mathbb{I}_{\Delta_k}(t).\] 
Given a metric space $U$, a $U$-valued random variable $x$ and a map $h:U\to \mathbb{R}^d$, the expectation of $h(x)$ is $\mathbb{E} h(x)$, and the covariance of $h(x)$ is $\mathbb{D} h(x)$. $L(t_0,T;U)$ denotes the space of $U$-valued Lebesgue measurable functions on $[t_0,T]$. 
We denote the set of positive semi-definite matrices with dimension $n$ as $S_{+}^{n}$. Given $A, B\in\mathbb{R}^{n\times n}$, we use $A \succeq B$ to indicate $A-B \in S_{+}^{n}$.

\textbf{Nonanticipative strategy with delay.}
Let $U$ and $V$ be the metric spaces in which player $1$ and player $2$ take control on the time interval $[t_0, T]$, resp. The admissible control spaces are $\mathcal{U}_p\triangleq L(t_0,T;U)$ and $\mathcal{V}_p\triangleq L(t_0,T;V)$ for each player. 
Classical nonanticipative strategies require players to act based solely on past information without forecasting future opponent behavior. However, physical and computational constraints can cause reaction lags, so the concept of NAD is introduced to model such buffer periods.
\begin{definition}[NAD {\cite[p.386]{basarHandbookDynamicGame2018}}]\label{def:NAD}
    A strategy $\alpha$ is called a nonanticipative strategy with delay for player $1$ if there exists a partition $\pi=\{t_0<\ldots<t_N=T\}$ such that, for any two control functions $v_1,v_2\in L(t_0,T;V)$ and any $k\in\{1,\ldots,N\}$, if $v_1= v_2$ a.e. on $[t_0,t_{k-1}]$, there is $\alpha(v_1)= \alpha(v_2)$ a.e. on $[t_0,t_{k}]$. 
    The sets of NADs for player $1$ and $2$ are denoted by $\mathcal{A}_{p}$ and $\mathcal{B}_{p}$, resp.
\end{definition}

 NAD strategies ensure a game is well-posed that, when one of the players adopts an NAD strategy, the interdependence of the control choices is resolved uniquely.
\begin{proposition}[Normal form of strategy {\cite[p.386]{basarHandbookDynamicGame2018}}]\label{prop:normal-form}
    Let $\alpha, \beta$ be nonanticipative strategies for player $1$ and $2$, resp. If at least one of these strategies is NAD, then there exists a unique pair of control functions $(u,v)\in L(t_0,T;U)\times L(t_0,T;V)$ such that $u = \alpha(v)$ and $v =\beta(u)$, a.e. on $[t_0,T]$.
\end{proposition}

\textbf{Weak approximation.}
Given a partition $\pi$ on a time horizon $[t_0,T]$, weak approximation characterizes the distributional deviation between a continuous-time stochastic process $\{\tilde{x}(t)\}_{t\in[t_0,T]}$ and a discrete-time stochastic process $\{x(k)\}_{k=0}^{|\pi|}$.
It requires that the probability distributions of $x(k)$ and $\tilde{x}(t_k)$ can be uniformly close for  $k\in\{0,1,\ldots,|\pi|\}$. Specifically, the difference between $\mathbb{E}\psi(\tilde{x}({t_k}))$ and $\mathbb{E}\psi(x(k))$ is small for any polynomial-growth function $\psi$. To make this statement rigorous, we first define the polynomial growth function.
\begin{definition}[Polynomial growth function]
	Let $\operatorname{PG}$ denote the set of continuous functions $\mathbb{R}^d \rightarrow \mathbb{R}$ of at most polynomial growth, i.e., $\psi(\cdot) \in \operatorname{PG}$ if there exists positive integers $\kappa_1, \kappa_2>0$ such that \[|\psi(x)| \leq \kappa_1\left(1+|x|^{\kappa_2}\right)\] for all $x \in \mathbb{R}^d$. 
     In particular, when $\kappa_2 = 1$, the above inequality is called the linear growth condition. For each integer $n \geq 1$ we denote by $\operatorname{PG}^n$ the set of $n$-times continuously differentiable functions $\mathbb{R}^d \rightarrow \mathbb{R}$ which, together with its partial derivatives up to and including order $n$, belong to $\operatorname{PG}$.
\end{definition}
It follows the complete definition of weak approximation.
\begin{definition}[Weak approximation]
	Given a partition $\pi$ on the time horizon $[t_0,T]$, a continuous-time stochastic process $\left\{\tilde{x}(t)\right\}_{t \in[t_0, T]}$ is an order $n$ weak approximation of a discrete-time stochastic process $\left\{x(k)\right\}_{k=0}^{|\pi|}$ if, $\forall \psi \in \operatorname{PG}^{n+1}$, there is a positive constant $C$ that does not depend on $\bar\pi$ such that
	\begin{equation*}
		\max _{k=0, \ldots, |\pi|}\left|\mathbb{E} \psi(x(k))-\mathbb{E} \psi(\tilde{x}(t_k))\right| \leq C \bar{\pi}^n.
	\end{equation*}
\end{definition}
\subsection{Zero-Sum Differential Game}
\textbf{Dynamics and cost.} Consider a ZSDG where both players' control processes are possibly influenced by a probability space $(\Omega,\mathcal{F},P)$. If neither of the control processes depends on the $\sigma$-field $\mathcal{F}$, the game is played under \textit{pure strategies}. Otherwise, the game is played under \textit{mixed strategies}. The dynamics and cost of the game are as follows,
\begin{equation*}\label{game:G}
	(\mathbf{G}):\left\{\begin{aligned}
		 & \dot{x}(t) = f(t,x(t), u(t), v(t)), \;  x(t_0) = x_0,\\
		 & J(t_0,x_0,u,v) = \mathbb{E}\int_{t_0}^{T}h(t,x,u,v)dt + g(x(T)).                            
	\end{aligned}\right.
\end{equation*}
where $x = \{x(t)\}_{t\in[t_0,T]}$ is the game state process taking values in $\mathbb{R}^d$. Given compact sets $U\subseteq\mathbb{R}^{d_1}$ and $V\subseteq\mathbb{R}^{d_2}$, the control processes of player $1$ and $2$ are $u:[t_0,T]\times\Omega\to U$ and $v:[t_0,T]\times\Omega\to V$, resp. 
For player $1$, the cost to be minimized is $J(t_0,x_0,u,v)$, and for player $2$, the cost to be minimized is $-J(t_0,x_0,u,v)$. 
The state evolution function $f:[t_0, T] \times \mathbb{R}^d \times U \times V \to \mathbb{R}^d$ is Lipschitz continuous and bounded. The terminal cost function $g: \mathbb{R}^{d} \to \mathbb{R}$ and the running cost function $h: [t_0,T]\times \mathbb{R}^d \times U \times V \to \mathbb{R}$ are Lipschitz continuous and bounded.
For an individual state trajectory under a realized sample, we define the accumulated cost function at time $t\in[t_0,T]$ as
\begin{equation}\label{eq:accumulated_cost}
    \operatorname{cost}(t)=\left\{\begin{array}{cl}
       \int_{t_0}^t h(s,x,u,v)ds  & t < T, \\ 
        \int_{t_0}^T h(s,x,u,v)ds + g(x(T)) & t=T,
    \end{array}\right.
\end{equation}
and there is $\mathbb{E}[\operatorname{cost}(T)] = J(t_0,x_0,u,v)$.

\textbf{Information and commitment pattern.} A strategy defines how a player selects admissible controls based on available information. The information pattern governs how players perceive the game, while the commitment pattern dictates how they respond.
\begin{assumption}[Information]\label{ass:info}
    During the game, the information available to both players includes the game dynamics, cost function, and history control process of its opponent.
\end{assumption}
Assumption \ref{ass:info} establishes the standard perfect information structure widely adopted in classical dynamic game literature \cite{basarDynamicNoncooperativeGame1998}. It ensures that players have access to the necessary system states and opponent histories to formulate rational feedback strategies.

\begin{assumption}[Commitment]\label{ass:commit}
    There is a commitment delay pattern $\pi = \{t_{0} < t_{1} < \ldots < t_{N} = T\}$, such that for any $k\in\{1,\ldots,N\}$, the actions of both players at time interval $[t_{k-1},t_k)$ are committed based on the information before $t_{k-1}$.
\end{assumption}
This assumption physically reflects the zero-order hold mechanism inherent in modern digital control systems caused by latencies. Furthermore, incorporating such a commitment delay is a standard mathematical paradigm for defining nonanticipative strategies in differential games \cite{buckdahnValueFunctionDifferential2013}.

Subject to the information and commitment patterns above, the strategy of player $1$ is a map from its opponent's admissible control space $\mathcal{V}$ to its own $\mathcal{U}$. Let $\mathcal{A}$ and $\mathcal{B}$ be the strategy spaces for player $1$ and $2$, resp. Given any pair of admissible strategies $(\alpha,\beta)\in(\mathcal{A},\mathcal{B})$, suppose that normal form exists, i.e., there is a unique couple of controls $(u,v)$ such that $\alpha(v)=u$ and $\beta(u)=v$. Then, one can define the associated cost functional $J(t_0,x_0,\alpha,\beta)$ as
\begin{equation}\label{eq:cost_functional}
    J(t_0,x_0,\alpha,\beta)\triangleq J(t_0,x_0,u,v).
\end{equation}

\textbf{Value functions.} To study the saddle-point equilibrium of a zero-sum differential game, the notions of lower and upper value functions are introduced.
\begin{definition}[Value functions]\label{def:value-function}
	Given admissible strategy spaces $\mathcal{A}$ and $\mathcal{B}$ for players $1$ and $2$, resp., the upper value function is
	\begin{equation}
		V_{+}(t_0, x_0)\triangleq\inf _{\alpha \in \mathcal{A}} \sup _{\beta \in \mathcal{B}} J\left(t_0, x_0, \alpha, \beta\right),
	\end{equation}
	and the lower value function is
	\begin{equation}
		V_{-}(t_0, x_0)\triangleq\sup_{\beta \in \mathcal{B}} \inf_{\alpha \in \mathcal{A}} J\left(t_0, x_0, \alpha, \beta\right).
	\end{equation}
\end{definition}
When the upper and lower value functions coincide, we say the saddle-point equilibrium exists under strategy space $\mathcal{A}$ and $\mathcal{B}$, and the game value function $V(t_0, x_0)$ is defined as either the upper or lower value function.
\begin{figure*}[t]
    \centering
    \includegraphics[width=0.88\linewidth]{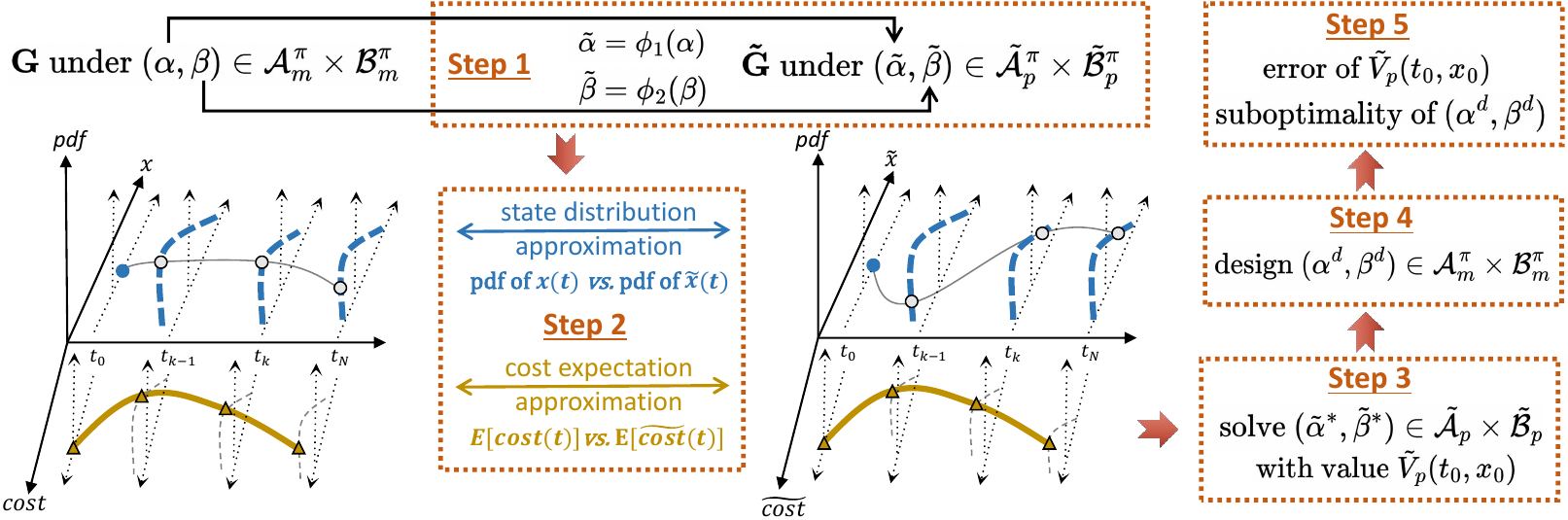}
    \caption{An illustration of five-step main methodology of the weak approximation framework.}
    \label{fig:methodology}
\end{figure*}
\subsection{Problems of Interest}
\textbf{Mixed strategy definition and existence of game value.} 
Existing definitions of mixed strategy for ZSDGs mainly focus on guaranteeing the existence of game value, where commitment delays are required to approach zero asymptotically.
For a ZSDG without the requirement of vanishing commitment delays, \textit{one needs a mixed strategy definition that guarantees the existence of game value and makes it practically possible to solve optimal mixed strategies}.

\textbf{Game dynamics and accumulated costs approximation.}
For a ZSDG played under mixed strategies, an immediate question arises: \textit{given a pair of admissible mixed strategies, how to characterize the dynamics of game state $x(t)$ and the accumulated cost $\operatorname{cost}(t)$?}  

\textbf{Near-optimal mixed strategy solution.}
Facilitated with a better characterization of the dynamics of a ZSDG under mixed strategies, our third question is \textit{how can we leverage this new understanding to solve optimal mixed strategies?} 

\section{Admissible Mixed Strategy: Definition, Properties, and Existence of Game Value}
This section rigorously defines the admissible mixed strategy spaces over non-vanishing commitment delays. We then explore their equivalent extensive forms to facilitate practical implementation, and establish the existence of the game value.

\subsection{Definition}
We first establish the underlying probability space $(\Omega,\mathcal{F},P)$ that governs the randomized behavior of both players:
\begin{equation}\label{eq:probability-space}
    (\Omega, \mathcal{F}, P):=\left(\left([0,1]^2\right)^{N}, \mathcal{B}\left([0,1]^2\right)^{\otimes N}, \lambda_2^{\otimes N}\right).
\end{equation}
Specifically, $\Omega=\left\{\left(\omega_{1,k}, \omega_{2,k}\right)_{k = 1}^{N} \mid \omega_{i,k} \in[0,1], i = 1,2 \right\}$ is a sample sapce containing all $[0,1]^2$-valued sequences. For player $i\in\{1,2\}$, the component $\omega_{i,k}$ represents the random seed influencing control actions during $[t_{k-1}, t_k)$. The probability measure $P=\lambda_2^{\otimes N}$ is constructed using the standard Lebesgue measure $\lambda_2$. Here, $P$ serves strictly as a fixed, foundational source of uniform randomness, rather than the strategic distribution itself.

To isolate the random-generation mechanism for player $i$ at step $k$, we define the local $\sigma$-field $\mathcal{G}_{k}^{(i)}\triangleq\sigma\{\omega_{i,k}\}$. The total history information up to step $k$ is captured by the filtration $\mathbb{F}\triangleq\{\mathcal{F}_k\}_{k=0}^\infty$, where $\mathcal{F}_k\triangleq\sigma\left\{\cup_{j \leq k}\left(\mathcal{G}_{j}^{(1)} \cup \mathcal{G}_{j}^{(2)}\right)\right\}$.

\begin{assumption}[Independence]\label{ass:independence}
    To ensure causality and fairness, the $\sigma$-fields $\mathcal{G}_{k}^{(1)}, \mathcal{G}_{k}^{(2)}$ and the history $\mathcal{F}_{k-1}$ are mutually independent for each $k \geq 1$.
\end{assumption}
Assumption \ref{ass:independence} requires that the local random noise generator $\mathcal{G}_k^{(i)}$ (e.g., the coin flip/seed) is independent of the past information $\mathcal{F}_{k-1}$. This ensures fairness, meaning neither player can magically foresee the opponent's private random seed.

\textbf{Admissible pure strategies.}
Given a metric space $U$, let $\operatorname{ST}(t_0,T;\pi,U)$ denote the space of $U$-valued stepwise functions governed by partition $\pi$, i.e,
\begin{equation}\label{eq:stepwise}
    \operatorname{ST}(t_0,T;\pi, U) \triangleq \left\{ u: u(t) \!=\!\sum\limits_{k=1}^{N}a_{k}\,\mathbb{I}_{\Delta_{k}}(t), a_{k} \in U \right\}.
\end{equation}
The admissible pure control space for player 1 is $\mathcal{U}_p^\pi\triangleq \operatorname{ST}(t_0,T;\pi,U)$, and the admissible pure strategy space is $\mathcal{A}_p^\pi\triangleq \{\alpha:\mathcal{V}_p^\pi\to\mathcal{U}_p^\pi\mid \alpha \text{ is NAD on }\pi\}$. Spaces $\mathcal{V}_p^\pi$ and $\mathcal{B}_p^\pi$ for player 2 are symmetrically defined. Proposition \ref{prop:normal-form} ensures a well-posed cost functional \eqref{eq:cost_functional}.

\textbf{Admissible mixed strategies.}
Let $L_{\mathbb{F}}(t_0,T;U)$ be the set of $\mathbb{F}$-adapted $U$-valued stochastic processes, and $\operatorname{ST}_{\mathbb{F}}(t_0,T;\pi,U)$ be its stepwise subspace.
\begin{definition}[Admissible control]\label{def:admissible-control}
    For player $1$, the admissible control space under mixed strategies is defined as 
    \begin{equation}
        \begin{aligned}
            \mathcal{U}_m^\pi \triangleq &\Big\{u \in \operatorname{ST}_{\mathbb{F}}(t_0,T;\pi,U) \mid \forall k\in\{1,\ldots,N\}, \\
            & u(t) \text{ is } \mathcal{F}_{k-1}\otimes\mathcal{G}_{k}^{(1)}\text{-measurable on } [t_{k-1},t_k)\Big\}.
        \end{aligned}
    \end{equation} 
    The space $\mathcal{V}^{\pi}_m$ for player 2 is similarly defined.
\end{definition}
It is crucial to distinguish the independent random seeds from the actual control actions. Assumption \ref{ass:independence} only enforces that the local noise generators (the "coin flips" $\mathcal{G}_{k}^{(i)}$) are independent of historical information, ensuring neither player can predict the opponent's private randomization. However, the generated control actions themselves are strictly dependent on the opponent's history. As formally defined in Definition \ref{def:admissible-control}, a player's admissible action is $\mathcal{F}_{k-1} \otimes \mathcal{G}_{k}^{(1)}$-measurable. It guarantees that a player leverages both their independent noise seed and the historical information $\mathcal{F}_{k-1}$, enabling them to dynamically react to the opponent's past decisions.

\begin{definition}[Admissible mixed strategy]\label{def:admissible-mixed-strategy}
    For player $1$, the admissible mixed strategy space is
    \begin{equation}\label{eq:mixed-strategy-map}
        \mathcal{A}_m^\pi \triangleq \{\alpha: \mathcal{V}^{\pi}_m \to \mathcal{U}^{\pi}_m\mid \text{subject to commitment delay } \pi\},
    \end{equation}
    where the delay $\pi$ enforces that for all $\mathbb{F}$-stopping times $\tau$ and controls $v_1,v_2\in\mathcal{V}^{\pi}_m$, if $v_1=v_2$ a.e. on $[t_0,t_{\tau}]\times\Omega$, then $\alpha(v_1)=\alpha(v_2)$ a.e. on $[t_0,t_{\tau+1}]\times\Omega$. $\mathcal{B}_m^{\pi}$ is similarly defined.
\end{definition}

For a ZSDG $\mathbf{G}$, we denote its upper and lower value functions under pure (resp. mixed) strategy spaces as $V^\pi_{p,+}, V^\pi_{p,-}$ (resp. $V^\pi_{m,+}, V^\pi_{m,-}$). Similar to \cite[Lemma 2.2]{buckdahnValueFunctionDifferential2013}, the normal form exists almost surely, yielding a well-defined expected cost $J(t_0,x_0,\alpha,\beta)$.

\subsection{Properties}
Directly endowing a functional space like $\mathcal{A}_p^\pi$ with a probability measure is notoriously difficult due to measure-theoretic obstacles \cite{aumannBorelStructuresFunction1961}. Instead, one can bridge pure and mixed strategies via a sample-path correspondence.
\begin{proposition}\label{prop:pure-and-mixed}
    Given $\omega\in\Omega$ and a mixed strategy $\alpha\in\mathcal{A}_m^\pi$, the corresponding pure strategy $\alpha^\omega$ defined via $\alpha^\omega(v(\cdot,\omega)) = \alpha(v)(\cdot,\omega)$ strictly satisfies $\alpha^\omega\in\mathcal{A}_p^\pi$.
\end{proposition}

While the normal form treats the strategy as a holistic map required for theoretical equilibrium proofs, practical execution necessitates a causal interpretation. To facilitate this, we introduce the extensive form, which unrolls the global strategy into a sequence of conditional random generators.
\begin{proposition}[Extensive form of mixed strategy]
    For player $1$, the extensive form of an admissible mixed strategy $\alpha\in\mathcal{A}_m^\pi$ is a sequence of $\mathcal{F}_{k-1}\otimes\mathcal{G}_{k}^{(1)}$-measurable maps $\{m_{k}^{(1)}\}_{k=1}^{|\pi|}$:
    \begin{equation}\label{def:mixed-extensive-form}
        \begin{aligned}
            m_{k}^{(1)}: \operatorname{ST}(t_{0},t_{k-1};\pi, V)\times\Omega &\to U\\
            \left(v_{\mid_{[t_0,t_{k-1})}},\omega\right) &\mapsto u(t_{k-1}).
        \end{aligned}	
    \end{equation}
\end{proposition} 
This confirms that, conditioned on historical information before $t_{k-1}$, players effectively sample static control actions for the upcoming interval using their local randomness $\mathcal{G}_k^{(i)}$. Consequently, the actual probability distribution of the mixed strategy over the control spaces is mathematically realized as the pushforward (induced) probability measure generated by the extensive form mapping $m_k^{(i)}$ under the fixed underlying measure $P$.

\subsection{Existence of Game Value} 
\begin{theorem}[Existence of game value]\label{thm:SE}
    The game $\mathbf{G}$ under the admissible mixed strategy spaces $\mathcal{A}_m^\pi\times\mathcal{B}_m^\pi$ has a unique game value, denoted as $V^\pi_m(t_0,x_0)$.
\end{theorem}
\begin{proof}
    Please see Appendix \ref{app:thm:SE}.
\end{proof}

\begin{remark}[Rationality of the strategy design]
    Our stepwise mixed strategy definition is both practically and theoretically motivated. Practically, it aligns perfectly with digital control execution (Zero-Order Hold) and circumvents normal-form paradoxes inherent to instantaneous nonanticipative strategies \cite[p.386]{basarHandbookDynamicGame2018}. Theoretically, as $\bar\pi \to 0$, stepwise functions densely approximate continuous $L_{\mathbb{F}}(t_0,T;U)$ spaces, allowing our framework to recover asymptotic mixed strategies \cite{buckdahnValueFunctionDifferential2013} while explicitly enabling the nonasymptotic approximation analysis developed in the following sections.
\end{remark}

\section{Weak Approximation Framework}
The key idea of the framework is straightforward: we design an SDG $\mathbf{\tilde{G}}$ under pure strategies that approximate $\mathbf{G}$ under mixed strategies with small errors in state distributions and cost expectations. Then we can solve the SDG $\mathbf{\tilde{G}}$ for game value and design near-optimal mixed strategies, with certified value approximation error and strategy suboptimality gap. Please see Fig. \ref{fig:methodology} for an illustration of the whole five steps that make up our weak approximation framework.

\subsection{Step 1: Design an SDG}
The goal of step 1 is to design a zero-sum SDG $\mathbf{\tilde{G}}$,
\begin{equation*}\label{game:tilde-G}
	(\mathbf{\tilde{G}}):\left\{\begin{aligned}
		 & d\tilde{x}(t) = \tilde{f}_1(t,\tilde{x}(t), \tilde{u}(t), \tilde{v}(t))dt \\
          & \quad+\tilde{f}_2(t,\tilde{x}(t), \tilde{u}(t), \tilde{v}(t))dW_t, \;  \tilde{x}(t_0) = \tilde{x}_0,\\
		 & \tilde{J}(t_0,\tilde{x}_0,\tilde{u},\tilde{v}) = \mathbb{E}\int_{t_0}^{T}\tilde{h}(t,\tilde{x},\tilde{u},\tilde{v})dt + \tilde{g}(\tilde{x}(T)),                    
	\end{aligned}\right.
\end{equation*}
where $\{W_t\}_{t\in[t_0,T]}$ is a standard $m$-dimensional Wiener process defined on a probability space $(\tilde{\Omega},\tilde{\mathcal{F}},\tilde{P})$, and the filtration is denoted as $\tilde{\mathbb{F}}\triangleq\{\tilde{\mathcal{F}}_t\}_{t\in[t_0,T]}$. $\tilde{U}$ and $\tilde{V}$ are compact metric spaces, the map $\tilde{f}_1:[t_0,T]\times\mathbb{R}^d\times\tilde{U}\times\tilde{V}\to\mathbb{R}^d$ and $\tilde{f}_2:[t_0,T]\times\mathbb{R}^d\times\tilde{U}\times\tilde{V}\to\mathbb{R}^{d\times m}$ are Lipschitz continuous and bounded. The cost functions $\tilde{g}$ and $\tilde{h}$ are Lipschitz continuous and bounded. Similarly to \eqref{eq:accumulated_cost}, one can define the accumulated cost function $\tilde{\operatorname{cost}}(t)$.
For player $1$, the admissible control space is
\begin{equation}\label{eq:tilde-admissible-control}
    \tilde{\mathcal{U}}_p^\pi\triangleq\operatorname{ST}_{\tilde{\mathbb{F}}}(t_0,T;\pi,\tilde{U}),
\end{equation}
and the admissible strategy space is
\begin{equation}\label{eq:tilde-admissible-strategy}
    \tilde{\mathcal{A}}_p^\pi\triangleq \{\tilde{\alpha}:\tilde{\mathcal{V}}_p^\pi\to\tilde{\mathcal{U}}_p^\pi\mid \tilde{\alpha} \text{ is NAD on }\pi\}.
\end{equation}
$\tilde{\mathcal{V}}_p$ and $\tilde{\mathcal{B}}_p^\pi$ for player $2$ are defined in a similar way.
 The upper and lower value functions of $\mathbf{\tilde{G}}$ under strategy spaces $\tilde{A}_p^\pi$ and $\tilde{B}_p^\pi$ are denoted as 
 \begin{equation}\label{eq:tilde-value-functions}
     \tilde{V}_{p,+}^\pi(t_0,x_0)\text{ and }\tilde{V}_{p,-}^\pi(t_0,x_0), \text{resp}.
 \end{equation}

It is important to highlight the theoretical connection between the information structures of $\mathbf{G}$ and $\mathbf{\tilde{G}}$ within this weak approximation framework. In standard stochastic dynamic games, strategies adapt to a common filtration. Our framework respects this paradigm: the pure strategies in the surrogate game $\mathbf{\tilde{G}}$ are fully adapted to the common filtration $\tilde{\mathbb{F}}$ generated by the stochastic state evolution. 
Furthermore, because Assumption \ref{ass:independence} guarantees that the original game prohibits arbitrary correlation between players' concurrent actions, the joint action distribution conditioned on the public history strictly factors into marginal distributions. Consequently, utilizing pure strategies in $\mathbf{\tilde{G}}$ adapted to the common filtration $\tilde{\mathbb{F}}$ is mathematically sufficient to fully characterize and approximate the conditional marginal distributions of the original mixed-strategy game, ensuring a rigorous structural alignment between the two games.
 
\subsection{Step 2: Certify SDG Weak Approximation of Order $n$}
The connection between $\mathbf{\tilde{G}}$ and $\mathbf{G}$ is made by designing a map from the admissible mixed strategy spaces of $\mathbf{G}$ to the admissible pure strategy spaces of $\mathbf{\tilde{G}}$. Then, we expect that $\mathbf{\tilde{G}}$ can approximate $\mathbf{G}$ in the sense of both state distributions and cost expectations with error-bound guarantees. 
\begin{definition}[SDG weak approximation]\label{def:SDG_weak_approximation}
    $\mathbf{\tilde{G}}$ under pure strategy spaces $\tilde{\mathcal{A}}_p^\pi\times\tilde{\mathcal{B}}_p^\pi$ is an order $n$ weak approximation of $\mathbf{G}$ under mixed strategy spaces $\mathcal{A}_m^\pi\times\mathcal{B}_m^\pi$ if, 
    \begin{enumerate}
        \item[(i)] there exist 
    \begin{equation}\label{eq:strategy-map}
        \phi_1: \mathcal{A}_m^\pi \to\tilde{\mathcal{A}}_p^\pi\text{ and }\phi_2: \mathcal{B}_m^\pi \to\tilde{\mathcal{B}}_p^\pi,
    \end{equation}
    where for any $(\alpha, \beta)\in\mathcal{A}_m^\pi\times\mathcal{B}_m^\pi$, one has $\tilde{\alpha}=\phi_1(\alpha)\in\tilde{\mathcal{A}}_p^\pi$ and $\tilde{\beta}=\phi_2(\beta)\in\tilde{\mathcal{B}}_p^\pi$. 
    \item[(ii)] the state trajectory of the SDG $\mathbf{\tilde{G}}$ under strategies $(\tilde{\alpha}, \tilde{\beta})$, i.e., $\{\tilde{x}(t)\}_{t\in[t_0,T]}$, is an order $n$ weak approximation of the state trajectory of the game $\mathbf{G}$ under strategies $(\alpha, \beta)$ on the time steps of $\pi$, i.e., $\{x(t_k)\}_{k=0}^{|\pi|}$.
    \item[(iii)] the accumulated cost trajectory of the SDG $\mathbf{\tilde{G}}$ under strategies $(\tilde{\alpha}, \tilde{\beta})$ approximates the accumulated cost trajectory of the game $\mathbf{G}$ under strategies $(\alpha, \beta)$ in an expected sense that $\left|\mathbb{E}\left[\operatorname{cost}(t_k)\right] - \mathbb{E}\left[\tilde{\operatorname{cost}}(t_k)\right]\right| = \mathcal{O}(\bar{\pi}^n)$ for $k=0,\ldots, |\pi|$.
    \end{enumerate}
\end{definition}

\subsection{Step 3: Solve $\mathbf{\tilde{G}}$} 
For a classical SDG, the admissible control spaces for both players are
\begin{equation}
    \tilde{\mathcal{U}}_p\triangleq\cup_{\pi\in\Pi_{[t_0,T]}}\,\tilde{\mathcal{U}}_p^\pi \text{ and } \tilde{\mathcal{V}}_p\triangleq\cup_{\pi\in\Pi_{[t_0,T]}}\tilde{\mathcal{V}}_p^\pi,
\end{equation}
and the admissible strategy spaces are
\begin{equation}
    \tilde{\mathcal{A}}_p \triangleq \cup_{\pi\in\Pi_{[t_0,T]}}\tilde{\mathcal{A}}_p^\pi \text{ and } \tilde{\mathcal{B}}_p \triangleq \cup_{\pi\in\Pi_{[t_0,T]}}\tilde{\mathcal{B}}_p^\pi.
\end{equation}
We denote the upper and lower value functions of $\mathbf{\tilde{G}}$ under strategy spaces $\tilde{\mathcal{A}}_p\times\tilde{\mathcal{B}}_p$ as
\begin{equation}
    \tilde{V}_{p,+}(t_0,x_0)\text{ and }\tilde{V}_{p,-}(t_0,x_0), \text{resp}.
\end{equation}
Usually, we assume that $\mathbf{\tilde{G}}$ is well-designed such that game value exists under $\tilde{\mathcal{A}}_p\times\tilde{\mathcal{B}}_p$. Therefore, one can solve the game value $\tilde{V}_p(t_0,x_0)$ and the optimal strategies $(\tilde{\alpha}^*,\tilde{\beta}^*)$.

\subsection{Step 4: Design Mixed Strategies}
Based on the optimal pure strategies for SDG $\mathbf{\tilde{G}}$ from step 3, the goal of step 4 is to design mixed strategies for $\mathbf{G}$. To this end, we first recall the definition of a zero-order holder for a stochastic process, i.e., $\operatorname{zoh}_{\pi}:L_{\mathbb{F}}(t_0,T;U)\to\operatorname{ST}_{\mathbb{F}}(t_0,T;\pi,U)$ satisfies that
\begin{equation*}
    \forall(t,\omega)\in[t_0,T]\times\Omega,\; \operatorname{zoh}_\pi[u](t,\omega)\triangleq \sum_{k=1}^N u(t_{k-1},\omega)\,\mathbb{I}_{\Delta_k}(t).
\end{equation*} 
Then, we introduce the concept of the zero-order holder of strategies as follows.
\begin{definition}[ZOH of strategies]\label{def:ZOH_of_strategy}
    Given $\pi\in\Pi_{[t_0,T]}$, the zero-order hold of strategies in $\tilde{\mathcal{A}}_p$ (resp., $\tilde{\mathcal{B}}_p$) is a map
    \begin{equation}
        \operatorname{ZOH}_\pi: \tilde{\mathcal{A}}_p \to \tilde{\mathcal{A}}_p, \quad \tilde{\alpha}\mapsto\operatorname{ZOH}_{\pi}[\tilde{\alpha}] 
    \end{equation}
    satisfies that $\forall\,\tilde{v}\in\tilde{\mathcal{V}}_p$ and $k\in\{1,\ldots,N\}$, there is 
    \begin{equation}
        \operatorname{ZOH}_{\pi}[\tilde{\alpha}](\tilde{v})(t,\omega)=\tilde{\alpha}(\tilde{v})(t_{k-1},\omega), 
    \end{equation}
    holds a.e. on $[t_{k-1},t_{k}]\times\Omega$.
\end{definition}
Notice that $\forall \tilde{\alpha}\in\tilde{\mathcal{A}}_p$, the strategy $\operatorname{ZOH}_{\pi}[\tilde{\alpha}]$ constrained to $\tilde{\mathcal{U}}_p^\pi$ belongs to $\tilde{\mathcal{A}}_p^\pi$, i.e., $\operatorname{ZOH}_{\pi}[\tilde{\alpha}]_{\mid_{\tilde{\mathcal{U}}_p^\pi}}\in\tilde{\mathcal{A}}_p^\pi$.
Therefore, one can apply the zero-order hold method to $(\tilde{\alpha}^*,\tilde{\beta}^*)$ and design near-optimal mixed strategies $(\alpha^d,\beta^d)\in\mathcal{A}_m^\pi\times\mathcal{B}_m^\pi$ satisfying 
\begin{equation}\label{eq:approximated-mixed-strategy-for-G}
    \left(\phi_1(\alpha^d),\phi_2(\beta^d)\right) = \left(\operatorname{ZOH}_\pi[\tilde{\alpha}^*]_{\mid_{\tilde{\mathcal{U}}_p^\pi}}, \operatorname{ZOH}_\pi[\tilde{\beta}^*]_{\mid_{\tilde{\mathcal{V}}_p^\pi}}\right).
\end{equation}

\subsection{Step 5: Performance Certification}
Finally, to rigorously certify the performance of the proposed methodology, we evaluate the bounds on the value approximation error and the suboptimality gap of the synthesized mixed strategies.

Suppose that $\mathbf{\tilde{G}}$ under $\tilde{\mathcal{A}}_p^\pi\times\tilde{\mathcal{B}}_p^\pi$ is an $n$-th order weak approximation of $\mathbf{G}$ under $\mathcal{A}_p^\pi\times\mathcal{B}_p^\pi$. Let $(\alpha^*,\beta^*)$ be optimal mixed strategies for $\mathbf{G}$, and $(\tilde{\alpha}^*, \tilde{\beta}^*)$ be optimal pure strategies for $\mathbf{\tilde{G}}$. Let $\tilde{\alpha}^{br} = \arg\min_{\tilde{\alpha}\in\tilde{\mathcal{A}}_p^{\pi}} \tilde{J}(t_0, x_0, \tilde{\alpha}, \operatorname{ZOH}_{\pi}[\tilde{\beta}^*])$ and $\tilde{\beta}^{br} = \arg\max_{\tilde{\beta}\in\tilde{\mathcal{B}}_p^{\pi}} \tilde{J}(t_0, x_0, \operatorname{ZOH}_{\pi}[\tilde{\alpha}^*], \tilde{\beta})$ denote the continuous-time best responses against the ZOH-executed opponent strategies in the surrogate game.

To facilitate our certification, we define the following error terms $\epsilon_1, \dots, \epsilon_5$:
\begin{equation*}
    \begin{aligned}
        \epsilon_1 =& |\tilde{J}(t_0,x_0,\tilde{\alpha}^*,\phi_2(\beta^*))-\tilde{J}(t_0,x_0,\operatorname{ZOH}_{\pi}[\tilde{\alpha}^*],\phi_2(\beta^*))|,\\
        \epsilon_2 =& |\tilde{J}(t_0,x_0,\phi_1(\alpha^*),\tilde{\beta}^*)-\tilde{J}(t_0,x_0,\phi_1(\alpha^*),\operatorname{ZOH}_{\pi}[\tilde{\beta}^*])|,\\
        \epsilon_3 =& |\tilde{J}(t_0,x_0,\tilde{\alpha}^*,\tilde{\beta}^*)-\tilde{J}(t_0,x_0,\operatorname{ZOH}_{\pi}[\tilde{\alpha}^*],\operatorname{ZOH}_{\pi}[\tilde{\beta}^*])|,\\
        \epsilon_4 =& |\tilde{J}(t_0, x_0, \tilde{\alpha}^{br}, \operatorname{ZOH}_{\pi}[\tilde{\beta}^*]) - \tilde{J}(t_0, x_0, \tilde{\alpha}^{br}, \tilde{\beta}^*)|,\\
        \epsilon_5 =& |\tilde{J}(t_0, x_0, \operatorname{ZOH}_{\pi}[\tilde{\alpha}^*], \tilde{\beta}^{br}) - \tilde{J}(t_0, x_0, \tilde{\alpha}^*, \tilde{\beta}^{br})|.
    \end{aligned}
\end{equation*}
Physically, these $\epsilon$ terms capture the non-smooth execution penalties caused by ``freezing'' the continuous-time optimal strategies into piecewise-constant ZOH actions. While analytically bounding them requires unrealistic global smoothness assumptions, they are practically negligible in modern high-frequency control systems. 

Furthermore, the maximum commitment delay $\bar{\pi}$ directly corresponds to the control sampling interval (e.g., $\bar{\pi}=0.05$ s for a $20\text{Hz}$ system), naturally ensuring $\bar{\pi} \ll 1$. As a result, the theoretical approximation error $\mathcal{O}(\bar{\pi}^n)$ shrinks rapidly with a higher weak approximation order $n$, yielding strictly tighter performance guarantees.

With these execution penalties isolated, we formally state the bounds on the value approximation error and the strategy suboptimality gaps.

\begin{theorem}[Value approximation]\label{thm:value_approximation}
    The game value $\tilde{V}_p(t_0,x_0)$ of $\mathbf{\tilde{G}}$ serves as a certified approximation to the true game value $V_m^\pi(t_0,x_0)$ of $\mathbf{G}$, with the error bounded by:
    \begin{equation}
        |V_m^\pi(t_0,x_0)-\tilde{V}_p(t_0,x_0)| \leq \max\{\epsilon_1, \epsilon_2\} + \mathcal{O}(\bar\pi^n).
    \end{equation}
\end{theorem}
\begin{proof}
    Please see Appendix \ref{app:thm:value_approximation}.
\end{proof}

\begin{theorem}[Suboptimality gap]\label{thm:suboptimality_gap}
    The near-optimal mixed strategies $(\alpha^d, \beta^d)$ possess the following bounded suboptimality gaps against the worst-case best responses:
    \begin{equation}
        \begin{aligned}
            J(t_0,x_0,\alpha^d,\beta^d) - \min_{\alpha\in\mathcal{A}_m^{\pi}}J(t_0, x_0, \alpha, \beta^d) \leq & \epsilon_3 + \epsilon_4 + \mathcal{O}(\bar{\pi}^n),\\
            \max_{\beta\in\mathcal{B}_m^{\pi}}J(t_0,x_0,\alpha^d,\beta)-J(t_0, x_0, \alpha^d, \beta^d) \leq& \epsilon_3 + \epsilon_5 + \mathcal{O}(\bar{\pi}^n).
        \end{aligned}
    \end{equation}
\end{theorem}
\begin{proof}
    Please see Appendix \ref{app:thm:suboptimality_gap}.
\end{proof}

\section{Control-Space Discretization Algorithm}
In this section, we apply the five-step theoretical framework proposed in Sec. V to general zero-sum differential games. While the weak approximation paradigm provides the necessary mathematical foundation for matching state distributions, finding optimal probability measures in continuous action spaces is generally intractable. To bridge the gap between abstract theory and algorithmic realization, we develop a control-space discretization algorithm. This approach maps the infinite-dimensional measure optimization problem into a sequence of solvable local linear programs, enabling the practical synthesis of near-optimal mixed strategies with certified performance.

\subsection{Step 1: Design an SDG $\mathbf{\tilde{G}}$}
The key idea of our control-space-discretization-based SDG design is to reduce the dimension of optimizing over continuous-state probability measures by optimizing over discrete-state probability measures.

\textbf{Discretized control space.}
Instead of optimizing over the infinite-dimensional probability measure spaces $\mathcal{P}(U)$ and $\mathcal{P}(V)$, we discretize the compact action spaces into finite sets of representative control points:
\begin{equation*}
    \begin{aligned}
        U_{\delta} =& \{u_1, u_2, \dots, u_{N_u}\} \subset U,\\
        V_{\delta} =& \{v_1, v_2, \dots, v_{N_v}\} \subset V.
    \end{aligned}
\end{equation*}
To explicitly quantify the resolution of this spatial approximation, we define the covering radii $\delta_u$ and $\delta_v$ for the discretized sets $U_{\delta}$ and $V_{\delta}$, respectively, such that:
\begin{equation*}
    \delta_u \triangleq \sup_{u \in U} \min_{u_i \in U_\delta} \|u - u_i\|, \quad \delta_v \triangleq \sup_{v \in V} \min_{v_j \in V_\delta} \|v - v_j\|.
\end{equation*}
Consequently, the continuous mixed strategies are parameterized as finite-dimensional probability mass functions (PMFs) evolving over time:
\begin{equation*}
    \begin{aligned}
        \tilde{u}(t) =& (\tilde{u}_1(t), \dots, \tilde{u}_{N_u}(t))^\top \in \Delta^{N_u-1}, \\
        \tilde{v}(t) =& (\tilde{v}_1(t), \dots, \tilde{v}_{N_v}(t))^\top \in \Delta^{N_v-1},
    \end{aligned}
\end{equation*}
where $\Delta^{N_u-1}$ denotes the standard $(N_u-1)$-dimensional probability simplex, strictly defined as:
\begin{equation*}
    \Delta^{N_u-1} \triangleq \left\{ p \in \mathbb{R}^{N_u} \;\middle|\; \sum_{i=1}^{N_u} p_i = 1, \; p_i \geq 0, \; \forall i \right\}.
\end{equation*}
This structural reduction is characterized by a natural parameterization mapping from the original infinite-dimensional measure space to the finite-dimensional simplex.

\textbf{Weakly approximated dynamics.}
Let $\{W_t\}_{t\in[t_0,T]}$ be a standard $d$-dimensional Wiener process defined on a probability space $(\tilde{\Omega},\tilde{\mathcal{F}},\tilde{P})$ with filtration $\tilde{\mathbb{F}}$. The continuous-time control spaces for $\mathbf{\tilde{G}}$ are defined as the standard simplices $\tilde{U}=\Delta^{N_u-1}$ and $\tilde{V}=\Delta^{N_v-1}$. Then, the pure strategy controls $\tilde{u}$ and $\tilde{v}$ are defined as the probability mass functions $\{(\tilde{u}_1(t), \dots, \tilde{u}_{N_u}(t))^\top\}_{t\in[t_0,T]} \in \tilde{U}$ and $\{(\tilde{v}_1(t), \dots, \tilde{v}_{N_v}(t))^\top\}_{t\in[t_0,T]} \in \tilde{V}$, resp.

Driven by $\tilde{u}$ and $\tilde{v}$, the zero-sum SDG $\mathbf{\tilde{G}}$ is governed by the following continuous-time stochastic differential equation:
\begin{equation*}
    d\tilde{x}(t) = \tilde{f}(t, \tilde{x}, \tilde{u}, \tilde{v})dt + \tilde{\Sigma}^{\frac{1}{2}}(t, \tilde{x}, \tilde{u}, \tilde{v}) dW_t, \quad \tilde{x}(t_0) = x_0,
\end{equation*}
where the equivalent drift vector $\tilde{f}$ and the diffusion matrix $\tilde{\Sigma}$ are synthesized by taking expectations over the discrete probability grids:
\begin{equation*}
    \begin{aligned}
        \tilde{f}(t, \tilde{x}, \tilde{u}, \tilde{v}) &= \sum_{i=1}^{N_u} \sum_{j=1}^{N_v} \tilde{u}_i(t) \tilde{v}_j(t) f(t, \tilde{x}, u_i, v_j), \\
        \tilde{\Sigma}(t, \tilde{x}, \tilde{u}, \tilde{v}) &= \sum_{i=1}^{N_u} \sum_{j=1}^{N_v} \tilde{u}_i(t) \tilde{v}_j(t) \Sigma(t, \tilde{x}, u_i, v_j).
    \end{aligned}
\end{equation*}
Here, $\Sigma(t, x, u, v) \in S_{+}^d$ acts as a \textit{designable} diffusion matrix function. We highlight a vital structural insight of this framework: while the drift term $\tilde{f}$ is strictly fixed to reflect the exact expectation of the system dynamics, the diffusion term $\tilde{\Sigma}$ provides a crucial degree of design flexibility. 

To mathematically guarantee that $\mathbf{\tilde{G}}$ remains a valid order $1$ weak approximation, the artificial noise injected into the system must scale with the commitment delay. Thus, as long as the underlying diffusion matrix $\Sigma$ is bounded and proportional to $\delta_\pi(t)$, the approximation order is theoretically preserved. This property will be formally proved and rigorously established in Step 2. 

A natural and mathematically robust choice for this design is the uncentered second moment, i.e., 
\begin{equation}\label{eq:simga_0}
    \Sigma_0(t, x, u, v) \triangleq \delta_\pi(t) f(t, x, u, v)f(t, x, u, v)^\top,
\end{equation}
which perfectly preserves Isaacs' condition for the resulting SDG $\mathbf{\tilde{G}}$. Nevertheless, other bounded functional designs for $\Sigma$ are also viable and can induce different empirical behaviors (e.g., varying degrees of conservatism or state exploration). We will comprehensively investigate and compare the practical performance of these alternative diffusion designs in the experimental section.

\textbf{Approximated cost functional.}
Correspondingly, the objective functional for $\mathbf{\tilde{G}}$ is reformulated to evaluate the expected cost under the discretized distributions:
\begin{equation*}
    \tilde{J}(t_0, x_0, \tilde{u}, \tilde{v}) = \mathbb{E} \left[ \int_{t_0}^T \tilde{h}(t, \tilde{x}, \tilde{u}, \tilde{v}) dt + \tilde{g}(\tilde{x}(T)) \right],
\end{equation*}
where the local running cost $\tilde{h}$ and the terminal cost $\tilde{g}$ are analytically given by:
\begin{equation*}
    \begin{aligned}
        \tilde{h}(t, \tilde{x}, \tilde{u}, \tilde{v}) &= \sum_{i=1}^{N_u} \sum_{j=1}^{N_v} \tilde{u}_i(t) \tilde{v}_j(t) h(t, \tilde{x}, u_i, v_j), \\
        \tilde{g}(\tilde{x}(T)) &= g(\tilde{x}(T)).
    \end{aligned}
\end{equation*}
By substituting the infinite-dimensional probability measures with finite-dimensional control vectors $(\tilde{u}, \tilde{v})$, the original game is successfully relaxed into an SDG that strictly preserves the saddle-point convexity properties.

To summarize, a zero-sum SDG is designed as follows,
\begin{equation*}
    (\mathbf{\tilde{G}})\!:\!\left\{\begin{aligned}
        &d\tilde{x}(t) \!=\! \tilde{f}(t, \tilde{x}, \tilde{u}, \tilde{v})dt \!+\! \tilde{\Sigma}^{\frac{1}{2}}(t, \tilde{x}, \tilde{u}, \tilde{v}) dW_t,  \tilde{x}(t_0) \!=\! x_0,\\
        & \tilde{J}(t_0,x_0,\tilde{u},\tilde{v}) = \mathbb{E} \left[ \int_{t_0}^T \tilde{h}(t, \tilde{x}, \tilde{u}, \tilde{v}) dt + \tilde{g}(\tilde{x}(T)) \right].
    \end{aligned}\right.
\end{equation*}

\subsection{Step 2: Certify SDG Weak Approximation of Order $1$}
To rigorously establish the weak approximation, we must specify the strategy projection maps $(\phi_1, \phi_2)$ required by \eqref{eq:strategy-map}. We achieve this by quantizing the conditional distributions of the extensive-form strategies across the disparate probability spaces.

First, let $\{U_1, U_2, \dots, U_{N_u}\}$ be a Borel partition of the compact action space $U$, generated by the representative discrete points $U_\delta$, such that $u \in U_i \implies \|u - u_i\| \leq \delta_u$ for all $i=1,\dots,N_u$. The partition for $V$, denoted as $\{V_1, \dots, V_{N_v}\}$, is symmetrically defined and bounded by $\delta_v$.

Next, paralleling the original game, we explicitly formulate the extensive form of the surrogate pure strategies.
\begin{proposition}
    For player $1$, an admissible pure strategy $\tilde{\alpha} \in \tilde{\mathcal{A}}_p^\pi$ in $\mathbf{\tilde{G}}$ has an equivalent extensive form as a sequence of measurable maps $\{\tilde{p}_k^{(1)}\}_{k=1}^{|\pi|}$ such that
    \begin{equation}\label{eq:tilde_extensive_form}
        \begin{aligned}
            \tilde{p}_k^{(1)}: \operatorname{ST}(t_0,t_{k-1};\pi,\tilde{V}) \times \tilde{\Omega} &\to \tilde{U} = \Delta^{N_u-1}\\
            (\tilde{v}_{\mid_{[t_0,t_{k-1})}}, \tilde{\omega}) &\mapsto \tilde{u}(t_{k-1}),
        \end{aligned}
    \end{equation}
    where $\tilde{p}_k^{(1)}$ is $\tilde{\mathcal{F}}_{t_{k-1}}$-measurable. The extensive form $\{\tilde{p}_k^{(2)}\}_{k=1}^{|\pi|}$ for player 2 is similarly defined.
\end{proposition}
Unlike the original mixed strategy, the surrogate map $\tilde{p}_k^{(1)}$ adapts to the intrinsic Wiener process via $\tilde{\Omega}$ but inherently lacks the private random seed $\mathcal{G}_k^{(1)}$. To bridge this gap, we define the exact projection mapping.
\begin{definition}[Strategy projection map]\label{def:general_map}
    The projection map $\phi_1: \mathcal{A}_m^\pi \to \tilde{\mathcal{A}}_p^\pi$ transforms a continuous mixed strategy $\alpha$ into a discrete pure strategy $\tilde{\alpha}$. 
    
    At any decision time $t_{k-1}$ and concurrent state $x \in \mathbb{R}^d$, the mapped strategy $\tilde{\alpha}$ outputs a probability mass function $\tilde{u}(t_{k-1}) \in \Delta^{N_u-1}$, defined exactly by the state-conditional distribution of the original random action $m_k^{(1)}$:
    \begin{equation}\label{eq:projection_map_u}
        \tilde{u}(t_{k-1}) \triangleq 
        \begin{bmatrix}
            P\left( m_k^{(1)} \in U_1 \mid x(t_{k-1}) = x \right) \\
            \vdots \\
            P\left( m_k^{(1)} \in U_{N_u} \mid x(t_{k-1}) = x \right)
        \end{bmatrix}.
    \end{equation}
    By structuring the PMF via state-conditional expectations, $\phi_1$ seamlessly synchronizes decisions across the disparate sample spaces $\Omega$ and $\tilde{\Omega}$. Symmetrically defining $\phi_2$ for Player 2 constitutes the complete feasible mapping.
\end{definition}
\begin{remark}[Interpretation of the mapping mechanism]
    This projection rigorously achieves two mathematical decoupling tasks:
    \begin{itemize}
        \item[(i)] \textbf{Randomization decoupling:} Taking the state-conditional expectation cleanly integrates out the independent private random seed $\mathcal{G}_{k}^{(1)}$, elegantly compressing the randomized strategy into a deterministic PMF vector.
        \item[(ii)] \textbf{Path decoupling:} Conditioning the mapping solely on the concurrent state $x$ is theoretically justified because the game value is purely Markovian (see Appendix \ref{app:thm:SE}). This fundamentally decouples the current optimal decisions from the disparate historical trajectories generated by $\Omega$ and $\tilde{\Omega}$.
    \end{itemize}   
\end{remark}

Under the projection map $(\phi_1, \phi_2)$, we can now establish the fundamental approximation guarantee. To ensure a tight coupling between the temporal and spatial discretizations, we require the spatial covering radii to scale linearly with the commitment delay, i.e., $\delta_u, \delta_v = \mathcal{O}(\bar{\pi})$. Additionally, the designed diffusion matrix must be bounded such that $\|\Sigma(t, x, u, v)\| \leq C_\Sigma \bar{\pi}$ for some positive constant $C_\Sigma$.

\begin{theorem}[$\mathbf{\tilde{G}}$ weak approximation]\label{thm:general_weak_approximation}
    Under the designed projection maps $(\phi_1, \phi_2)$, the surrogate SDG $\mathbf{\tilde{G}}$ using pure strategies is an order $1$ weak approximation (i.e., $\mathcal{O}(\bar\pi)$) of the original general game $\mathbf{G}$ using mixed strategies.
\end{theorem}
\begin{proof}
    Please see Appendix \ref{app:thm:general_weak_approximation}.
\end{proof}

\subsection{Step 3: Solve $\mathbf{\tilde{G}}$}
Having constructed the approximated stochastic differential game $\mathbf{\tilde{G}}$ under the finite-dimensional pure strategy spaces $\tilde{\mathcal{A}}_p^\pi \times \tilde{\mathcal{B}}_p^\pi$, the problem of finding the game value and optimal strategies is transformed into solving its associated Hamilton-Jacobi-Isaacs (HJI) equation. 

For the continuous-time SDG $\mathbf{\tilde{G}}$, we define the Hamiltonian $\mathcal{H}: [t_0, T] \times \mathbb{R}^d \times \mathbb{R}^d \times S^d \times \Delta^{N_u-1} \times \Delta^{N_v-1} \to \mathbb{R}$ as:
\begin{equation*}
    \begin{aligned}
        \mathcal{H}(t, x, p, M, \tilde{u}, \tilde{v}) \triangleq &\tilde{h}(t, x, \tilde{u}, \tilde{v}) + p^\top \tilde{f}(t, x, \tilde{u}, \tilde{v})\\
        &\qquad + \frac{1}{2}\operatorname{Tr}\left[\tilde{\Sigma}(t, x, \tilde{u}, \tilde{v}) M\right].
    \end{aligned}
\end{equation*}\\
Recalling the definitions from Step 1, the functions $\tilde{f}$, $\tilde{\Sigma}$, and $\tilde{h}$ are formulated as expectations over the discrete probability grids. Consequently, the Hamiltonian is perfectly bilinear with respect to the probability mass functions $\tilde{u}$ and $\tilde{v}$. It can be rewritten as a standard matrix game:
\begin{equation*}
    \mathcal{H}(t, x, p, M, \tilde{u}, \tilde{v}) = \tilde{u}^\top Q(t, x, p, M) \tilde{v},
\end{equation*}
where the $(i, j)$-th entry of the payoff matrix $Q \in \mathbb{R}^{N_u \times N_v}$ is explicitly given by evaluating the original system dynamics and costs at the discrete representative points $(u_i, v_j)$:
\begin{equation*}
    Q_{i,j} \!=\! h(t, x, u_i, v_j) + p^\top f(t, x, u_i, v_j) + \frac{1}{2}\! \operatorname{Tr}\! \Big[ \Sigma(t, x, u_i, v_j) M \Big].
\end{equation*}\\
Because $\mathcal{H}$ is a continuous bilinear form mapping from the compact and convex standard simplices $\Delta^{N_u-1} \times \Delta^{N_v-1}$, von Neumann's Minimax Theorem strictly guarantees that Isaacs' condition holds. That is, the order of optimization is interchangeable:
\begin{equation}\label{eq:isaacs_condition}
    \min_{\tilde{u} \in \Delta^{N_u-1}} \max_{\tilde{v} \in \Delta^{N_v-1}} \tilde{u}^\top Q \tilde{v} = \max_{\tilde{v} \in \Delta^{N_v-1}} \min_{\tilde{u} \in \Delta^{N_u-1}} \tilde{u}^\top Q \tilde{v}.
\end{equation}\\
Consequently, the game value of $\mathbf{\tilde{G}}$ is governed by the following second-order nonlinear parabolic partial differential equation, known as the HJI equation:
\begin{equation}\label{eq:hji_pde}
    \left\{\begin{aligned}
        &-\partial_t \tilde{V}(t,x) = \min_{\tilde{u} \in \Delta^{N_u-1}} \max_{\tilde{v} \in \Delta^{N_v-1}} \mathcal{H}(t, x, \nabla_x \tilde{V}, \nabla_{xx}^2 \tilde{V}, \tilde{u}, \tilde{v}),\\
        &\tilde{V}(T, x) = g(x),
    \end{aligned}\right.
\end{equation}
for $(t, x) \in [t_0, T) \times \mathbb{R}^d$.

To rigorously establish the existence and uniqueness of the solution to \eqref{eq:hji_pde}, we appeal to the standard viscosity solution framework for stochastic differential games.

\begin{theorem}[Viscosity solution and optimal strategies]\label{thm:viscosity_solution}
    For the surrogate SDG $\mathbf{\tilde{G}}$, the following properties hold:
    \begin{enumerate}
        \item[(i)] \textbf{Game value:} It has a unique game value $\tilde{V}_p(t, x)$, which is the unique viscosity solution to the HJI equation \eqref{eq:hji_pde}.
        \item[(ii)] \textbf{Optimal strategies:} At any time-state pair $(t, x)$, the optimal control inputs $(\tilde{u}^*, \tilde{v}^*)$ correspond to the saddle point of the local matrix game $Q$, which can be directly solved via linear programming.
    \end{enumerate}
\end{theorem}
\begin{proof}
    Because $\mathbf{\tilde{G}}$ is formulated as a classical stochastic differential game driven by a Wiener process, where the players independently select actions from the compact metric spaces $\tilde{U}$ and $\tilde{V}$, the existence of the value and its characterization as the unique viscosity solution follow directly from the dynamic programming principle and the comparison principle established in the foundational literature of stochastic differential games (e.g., Fleming and Souganidis \cite{flemingExistenceValueFunctions1989b}). 
\end{proof}

The proposed finite-dimensional HJI equation \eqref{eq:hji_pde} bridges actionable algorithms with theoretical limits. Notably, Buckdahn et al. \cite{buckdahnValueFunctionDifferential2013} proved that as the commitment delay vanishes ($|\pi| \to 0$), the value of a mixed-strategy ZSDG converges to the viscosity solution of a first-order HJI equation. 

\begin{remark}[Connection to asymptotic mixed-strategy PDEs]
Our formulation advances the foundational result in \cite{buckdahnValueFunctionDifferential2013} from two key perspectives:
\begin{enumerate}
    \item[(i)] \textbf{Asymptotic consistency:} Because our designed diffusion matrix scales with the delay ($\tilde{\Sigma} \propto \delta_\pi(t)$), the second-order term $\frac{1}{2}\operatorname{Tr}(\tilde{\Sigma} \nabla_{xx}^2 \tilde{V})$ uniformly vanishes as $\bar{\pi} \to 0$. Concurrently, as $\delta_u, \delta_v \to 0$, the simplex optimization recovers the full measure spaces. Thus, our second-order PDE flawlessly collapses into the classical first-order PDE derived in \cite{buckdahnValueFunctionDifferential2013}.
    \item[(ii)] \textbf{Finite-delay solvability:} While existing literature characterizes value functions under an idealized continuous limit, our formulation explicitly preserves the physical delay $\delta_\pi(t)$ via the diffusion term. This critical mechanism transforms the intractable measure-optimization into solvable local matrix games, enabling the direct computation of executable mixed strategies.
\end{enumerate}    
\end{remark}

\subsection{Step 4-5: Strategy Synthesis and Performance Certification}
With the theoretical HJI equation \eqref{eq:hji_pde} established, we now operationalize the solution into an executable mixed strategy for the original game $\mathbf{G}$. 

Algorithm \ref{alg:g_general} details this synthesis process. At each commitment epoch $t_{k-1}$, the players construct the local payoff matrix $Q$ using the spatial gradients of the value function $\tilde{V}$ and the chosen designable diffusion matrix $\Sigma$. By solving a standard Linear Program (LP) to find the minimax saddle point, they obtain their optimal PMFs $\tilde{u}^*$ and $\tilde{v}^*$. The actual control actions are then randomly drawn from the discrete representative sets $U_\delta$ and $V_\delta$ according to these optimal probabilities and held constant over the interval (ZOH).

Having designed the practical execution algorithm, we formally certify its global performance bounds. By substituting the order $1$ fundamental weak approximation established in Theorem \ref{thm:general_weak_approximation} into the general certification framework of Section V, we derive the explicit performance guarantees. 

It is crucial to acknowledge that in general nonlinear settings with discretized action spaces, the theoretical continuous-time optimal strategies of $\mathbf{\tilde{G}}$ may exhibit bang-bang type discontinuities (i.e., jumping between simplex vertices) across switching surfaces in the state space. Therefore, the execution penalties introduced by the Zero-Order Hold (ZOH) mechanism cannot be trivially globally bounded by $\mathcal{O}(\bar{\pi})$, and are thus explicitly preserved.

\begin{algorithm}[t]
\renewcommand{\algorithmicrequire}{\textbf{Input:}}
\renewcommand{\algorithmicensure}{\textbf{Output:}}
\caption{Near-optimal mixed strategies synthesis for $\mathbf{G}$} 
\label{alg:g_general}
\begin{algorithmic}[1]
\Require Game dynamics $f$, cost $h$, terminal cost $g$, continuous action spaces $U, V$, commitment delay pattern $\pi$, spatial grids $U_\delta, V_\delta$, bounded diffusion matrix design $\Sigma$, and the solved HJI value function $\tilde{V}$.
\State Initialize the accumulated cost $J \leftarrow 0$.
\For {$k=1,\ldots,|\pi|$}
    \State Observe current game state $x(t_{k-1})$.
    \State Evaluate gradients $\nabla_x \tilde{V}$ and $\nabla_{xx}^2 \tilde{V}$ at $(t_{k-1}, x(t_{k-1}))$.
    \State Construct local payoff matrix $Q \in \mathbb{R}^{N_u \times N_v}$ where:
    \begin{equation*}
        Q_{i,j} = h_k + (\nabla_x \tilde{V})^\top f_k + \frac{1}{2} \operatorname{Tr} \left( \Sigma_k \nabla_{xx}^2 \tilde{V} \right),
    \end{equation*}
    with 
    \begin{equation*}
        \begin{aligned}
            &f_k \triangleq f(t_{k-1}, x(t_{k-1}), u_i, v_j)\\
            &h_k \triangleq h(t_{k-1}, x(t_{k-1}), u_i, v_j)\\
            &\Sigma_k \triangleq \Sigma(t_{k-1}, x(t_{k-1}), u_i, v_j).
        \end{aligned}
    \end{equation*}
    \State Solve the local Minimax LP to obtain optimal PMFs:
    \begin{equation*}
        (\tilde{u}^*, \tilde{v}^*) = \arg\min_{\tilde{u} \in \Delta^{N_u-1}} \max_{\tilde{v} \in \Delta^{N_v-1}} \tilde{u}^\top Q \tilde{v}.
    \end{equation*}
    \State Randomly sample actions from the spatial grids: 
    \begin{equation*}
        u(t_{k-1}) \sim \tilde{u}^* \quad \text{and} \quad v(t_{k-1}) \sim \tilde{v}^*.
    \end{equation*}
    \For {$t\in[t_{k-1},t_k)$}
        \State Apply ZOH: $u(t) = u(t_{k-1})$ and $v(t) = v(t_{k-1})$.
        \State True game state evolves: $$\dot{x}(t) = f(t, x(t), u(t), v(t)).$$
    \EndFor
    \State Update cost: $J \leftarrow J + \int_{t_{k-1}}^{t_k} h(t, x(t), u(t), v(t)) dt$.
\EndFor
\State $J \leftarrow J + g(x(T))$.
\Ensure Sampled control sequence $\{u(t_{k-1}), v(t_{k-1})\}_{k=1}^{|\pi|}$ and total cost $J$.  
\end{algorithmic} 
\end{algorithm}

\begin{theorem}[Performance certification]\label{thm:general_performance}
    The mixed strategies synthesized via $\mathbf{\tilde{G}}$ and executed by Algorithm \ref{alg:g_general} possess the following bounds for the original game $\mathbf{G}$:
    \begin{enumerate}
        \item[(i)] Value approximation error: 
    \begin{equation*}
        |V_m^\pi(t_0,x_0)-\tilde{V}_p(t_0,x_0)| \leq \max\{\epsilon_1, \epsilon_2\} + \mathcal{O}(\bar\pi).
    \end{equation*}
        \item[(ii)] Strategy suboptimality gap: 
    \end{enumerate}
    \begin{equation*}
        \left\{\begin{aligned}
            &J(t_0,x_0,\alpha^d,\beta^d) - \min_{\alpha\in\mathcal{A}_m^{\pi}}\!J(t_0, x_0, \alpha, \beta^d) \leq \epsilon_3 + \epsilon_4 + \mathcal{O}(\bar{\pi}),\\
            &\max_{\beta\in\mathcal{B}_m^{\pi}}\!J(t_0,x_0,\alpha^d,\beta)-J(t_0, x_0, \alpha^d, \beta^d) \leq \epsilon_3 + \epsilon_5 + \mathcal{O}(\bar{\pi}).
        \end{aligned}\right.
    \end{equation*}
\end{theorem}

\begin{remark}[Nature of the ZOH execution penalties]
    The explicit retention of the error terms $\epsilon_1 \dots \epsilon_5$ highlights a fundamental mathematical characteristic of optimizing over discrete spatial domains:
    \begin{enumerate}
        \item[(i)] \textbf{Bang-bang discontinuities:} Because the local Hamiltonian is strictly bilinear with respect to the probability vectors $\tilde{u}$ and $\tilde{v}$, the local optimal solutions are inherently bang-bang. This inevitably leads to discontinuous strategy jumps (e.g., abruptly shifting between simplex vertices) as the system crosses switching surfaces in the state space.
        \item[(ii)] \textbf{The ZOH execution gap:} The $\epsilon$ terms explicitly isolate and quantify the unavoidable execution penalties generated when these theoretically discontinuous optimal strategies are ``frozen'' in time by the practical piecewise-constant ZOH delay pattern $\pi$.
        \item[(iii)] \textbf{Conditions for vanishing penalties:} These penalties reflect the physical reality of discrete execution rather than a systemic flaw. If the system dynamics inherently avoid switching surfaces, or if algorithmic smoothing techniques (such as entropy regularization) are applied to the local matrix games, these $\epsilon$ penalties will strictly diminish toward $\mathcal{O}(\bar{\pi})$.
    \end{enumerate} 
\end{remark}

\section{Numerical Experiments}

In this section, we validate the proposed control-space-discretization-based weak approximation framework through a highly nonlinear zero-sum differential game. Specifically, we evaluate the state and cost approximations under different diffusion designs, demonstrate the necessity of mixed strategies by analyzing the suboptimality gaps against pure-strategy bounds, and finally, verify our theoretically established $\mathcal{O}(\bar{\pi})$ error bounds through a rigorous sensitivity analysis on the commitment delay $\pi$.

\subsection{Experimental Setup}
Consider a scalar nonlinear zero-sum differential game $\mathbf{G}$ defined on the time horizon $[t_0, T] = [0, 1]$ with the initial state $x_0 = 1.0$. The continuous-time system dynamics are driven by both players as follows:
\begin{equation}
    f(t,x,u,v) = \sin(x) - 0.5x + 2u\cos(x/4) - 2v,
\end{equation}
where the compact action spaces are $u \in U = [-1, 1]$ and $v \in V = [-1, 1]$. The running cost penalizes state deviations, action efforts, and intense strategy collisions:
\begin{equation}
    h(t,x,u,v) = 3x^2 + 40\exp(-15(u - v)^2) + u^2 - v^2,
\end{equation}
and the terminal cost is given by $g(x(T)) = 2x^2$. 

To apply the proposed methodology, we uniformly discretize the continuous action spaces into $N_u = N_v = 41$ grid points, and discretize the state space over $x \in [-\pi, \pi]$ into $63$ grid points. For the numerical integration of the game, we set the simulation time step to $dt = 0.001$. The default baseline commitment delay is set to an equispaced partition with $\bar{\pi} = 0.05$. 

We evaluate three distinct functional designs for the artificial diffusion term $\tilde{\Sigma}$ in the approximated SDG $\mathbf{\tilde{G}}$, corresponding to varying degrees of local variance matching:
\begin{enumerate}
    \item \textbf{Zero diffusion}: $\sigma_1 = 0$;
    \item \textbf{Absolute drift}: $\sigma_2 = |f(t,x,u,v)|$;
    \item \textbf{Delay-scaled drift}: $\sigma_3 = \sqrt{\delta_\pi(t)} |f(t,x,u,v)|$, which strictly adheres to the $\mathcal{O}(\bar{\pi})$ condition in Theorem \ref{thm:general_weak_approximation}.
\end{enumerate}
The approximated game values and optimal PMFs are obtained by solving the discrete HJI equation via backward dynamic programming. All stochastic evaluations are averaged over $2,000$ independent Monte Carlo sample paths.

\subsection{Weak Approximation Performance}
We first evaluate how well the designed pure-strategy SDG $\mathbf{\tilde{G}}$ approximates the original game $\mathbf{G}$ under the designed mixed strategies (using the theoretically sound design $\sigma_3$). 
\begin{figure*}[h]
    \centering
    \subfigure[zero diffusion]{
        \includegraphics[width = 0.31\textwidth]{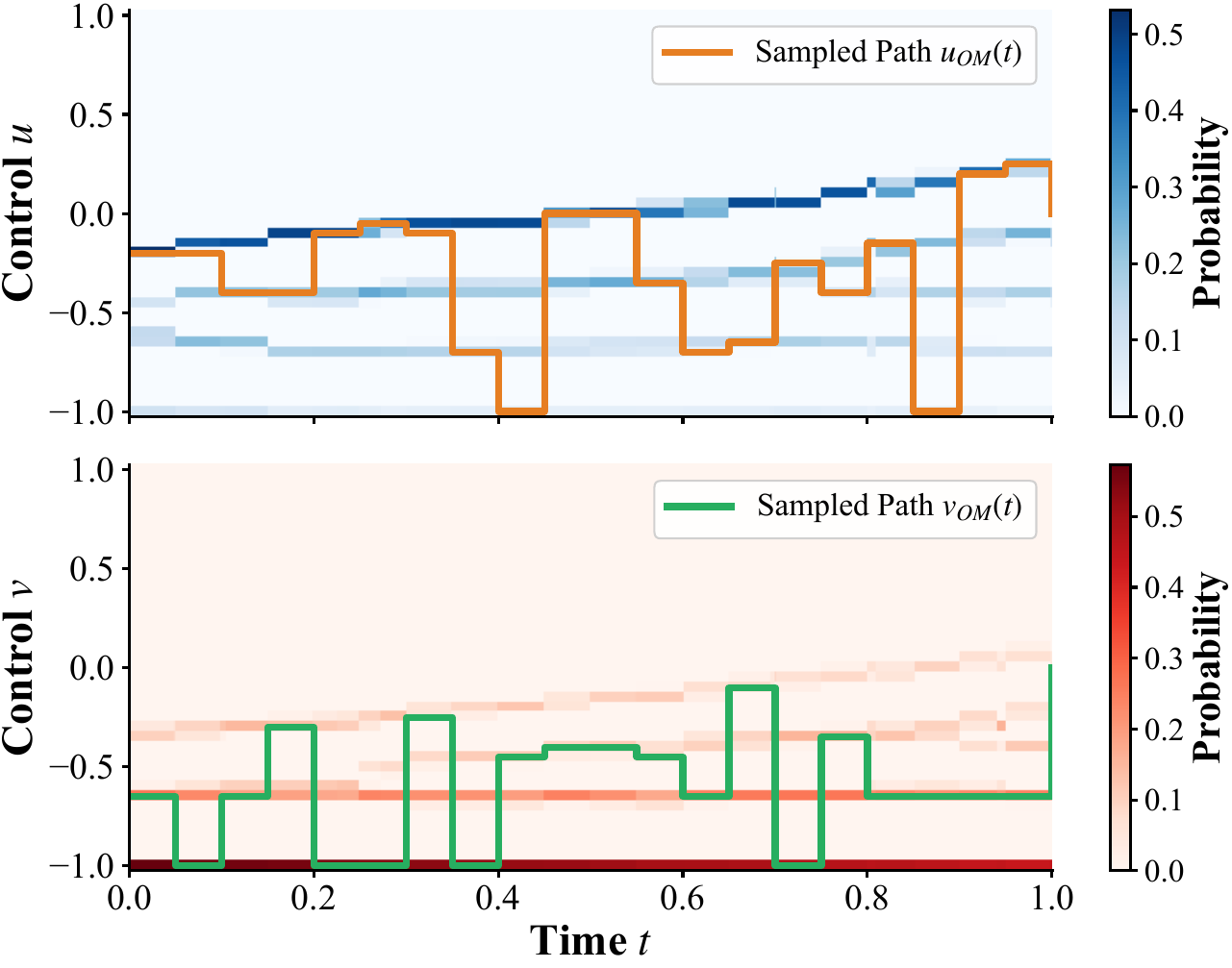}
    }
    \subfigure[absolute drift]{
        \includegraphics[width = 0.31\textwidth]{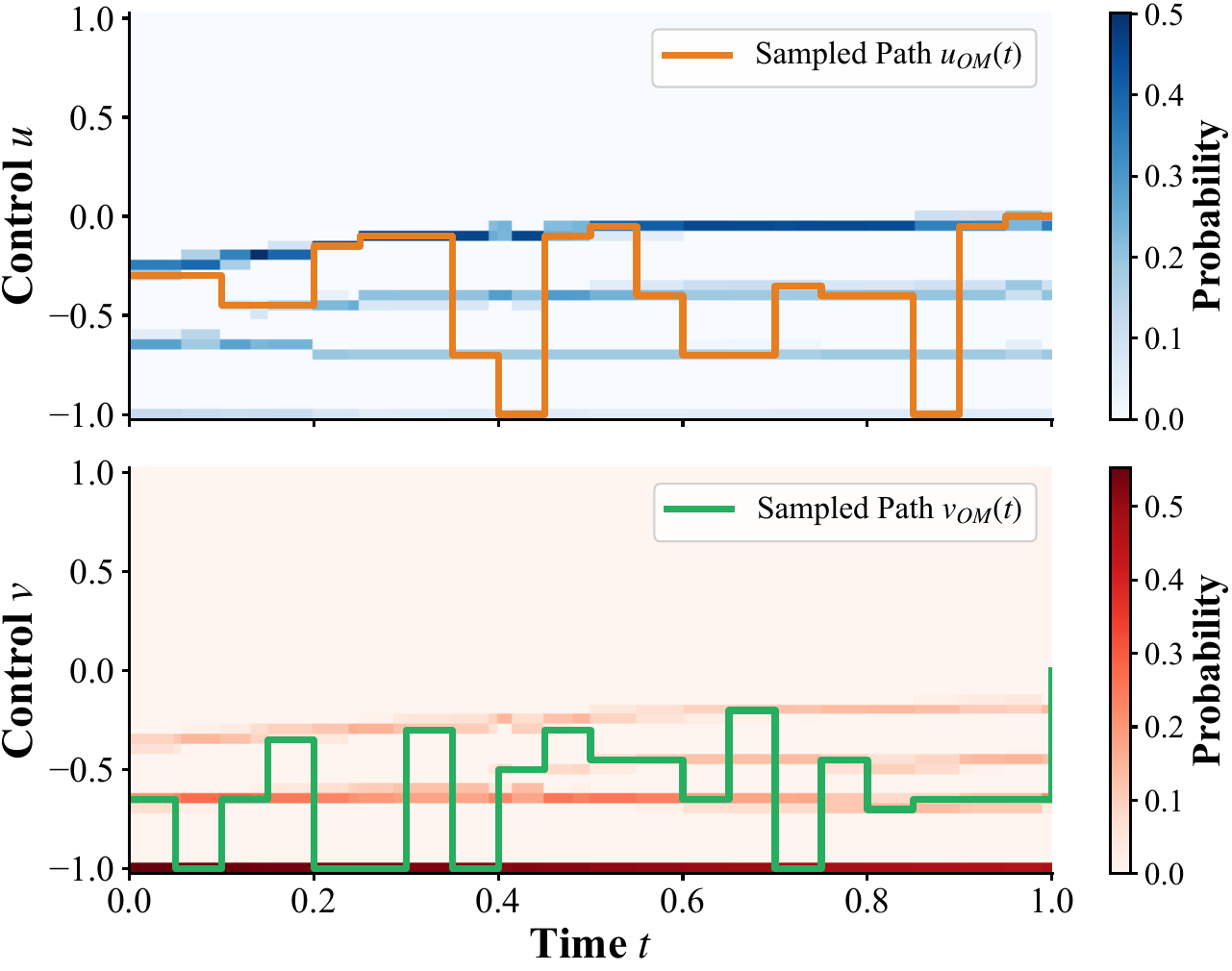}
    }
    \subfigure[delay-scaled drift]{
        \includegraphics[width = 0.31\textwidth]{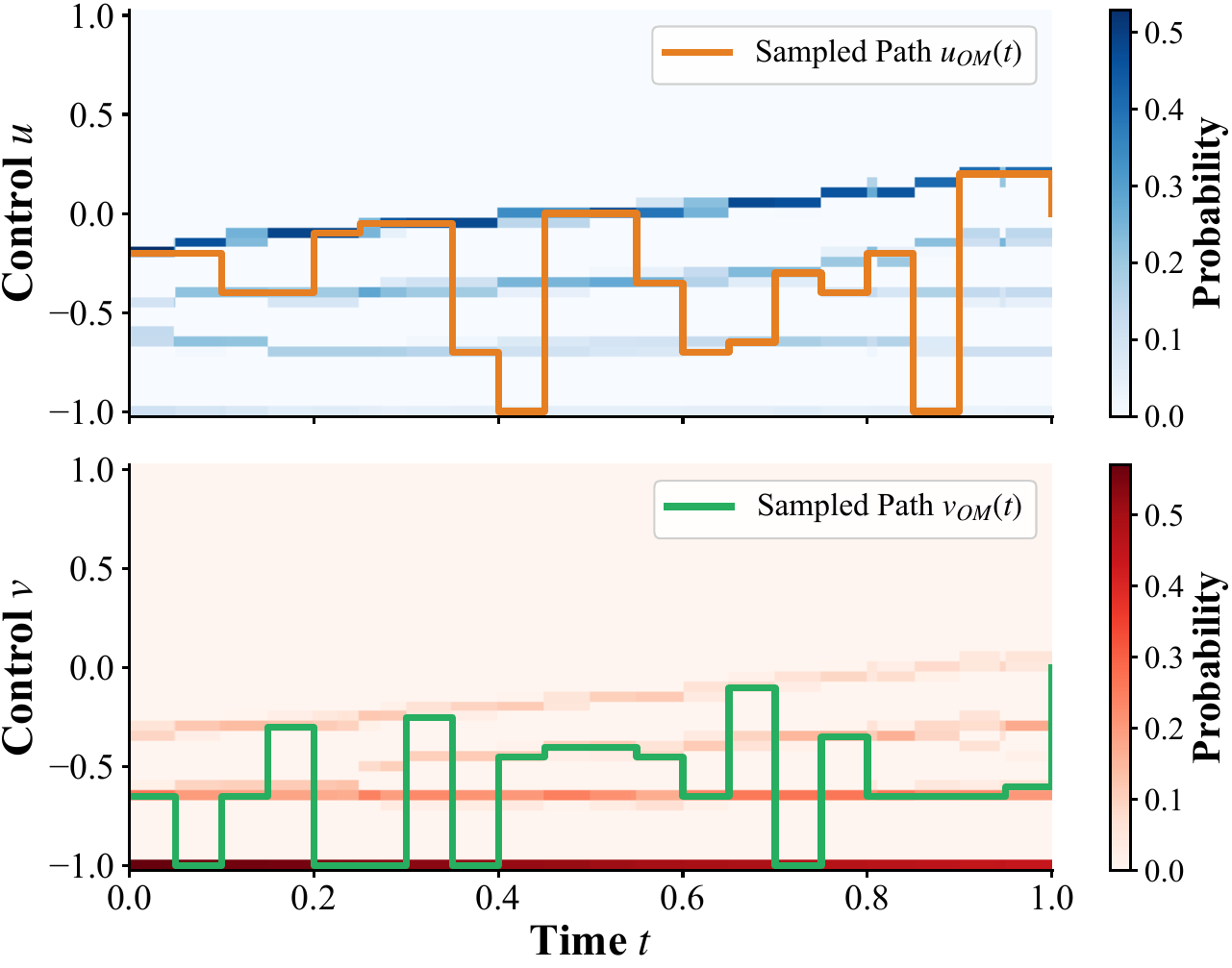}
    }
    \caption{Mixed strategy distribution under different designs for the diffusion term of $\mathbf{\tilde{G}}$.}
    \label{fig:heatmap}
\end{figure*}

Fig. \ref{fig:heatmap} visualizes the time-varying probability mass functions (PMFs) $\tilde{u}^*(t)$ and $\tilde{v}^*(t)$ solved from $\mathbf{\tilde{G}}$. The optimal mixed strategy exhibits a highly multimodal distribution to actively avoid the opponent's actions due to the severe collision penalty. The overlaid stepwise trajectory confirms that the Zero-Order Hold (ZOH) sampling effectively draws control commands governed by these optimal distributions.

\begin{figure*}[h]
    \centering
    \subfigure[weak approx: state Wasserstein distance]{
        \raisebox{0.3em}{\includegraphics[width = 0.3\textwidth]{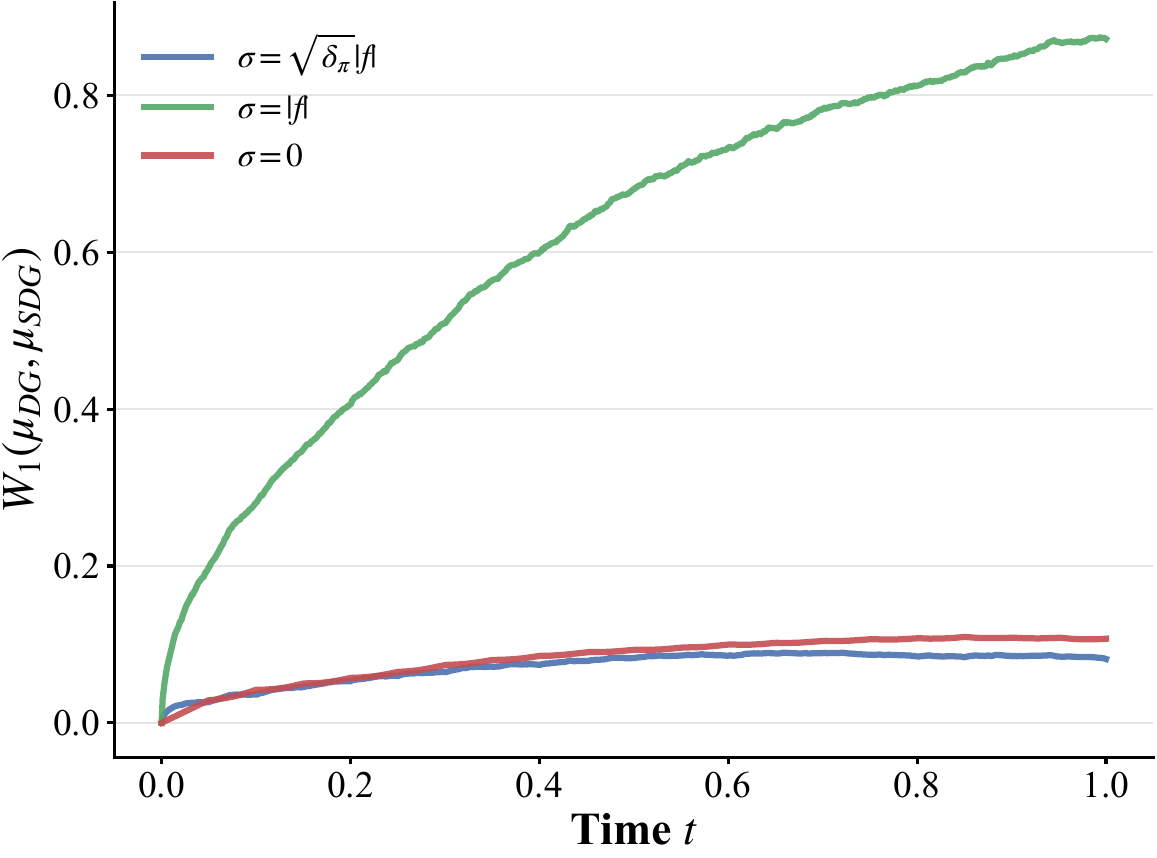}}
        \label{fig:wasserstein_compare}
    }
    \subfigure[weak approx: cost expectation alignment]{
        \includegraphics[width = 0.31\textwidth]{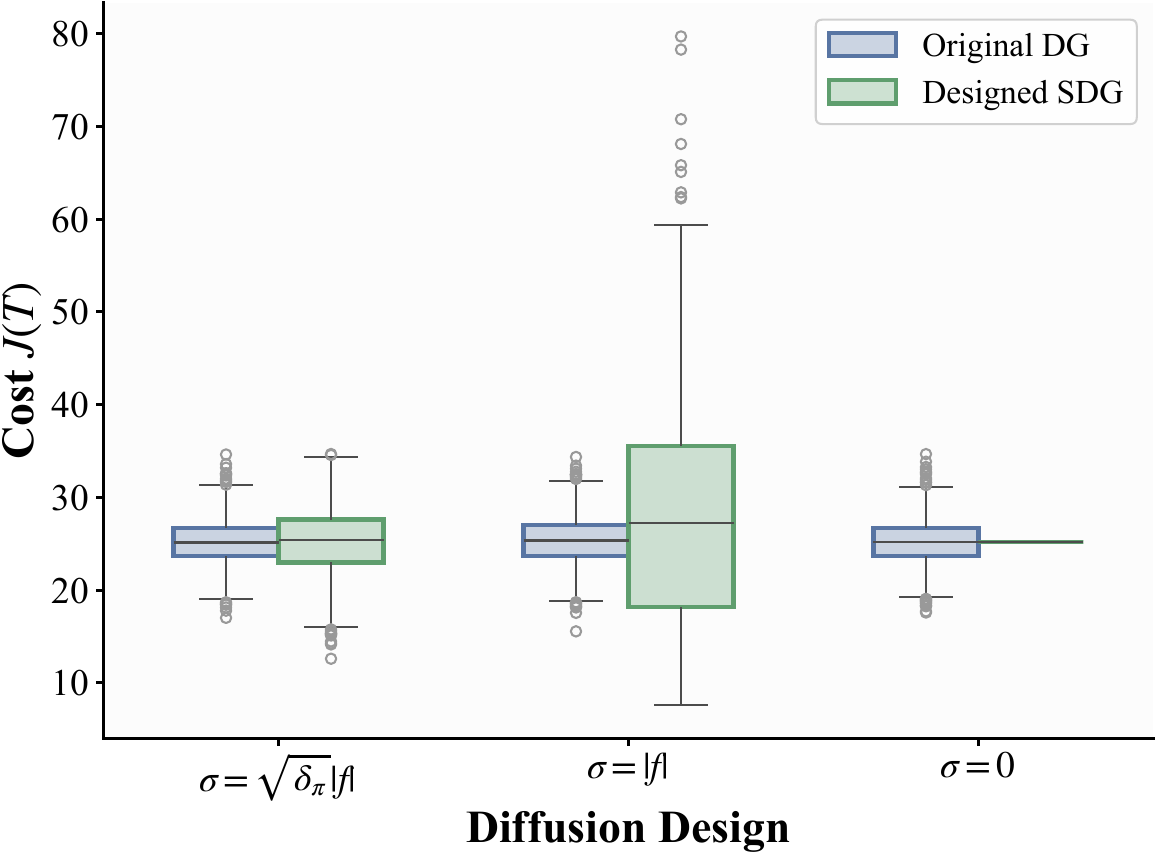}
        \label{fig:cost_boxplot}
    }
    \subfigure[suboptimality \& pure strategy value gap]{
        \raisebox{0.4em}{\includegraphics[width = 0.3\textwidth]{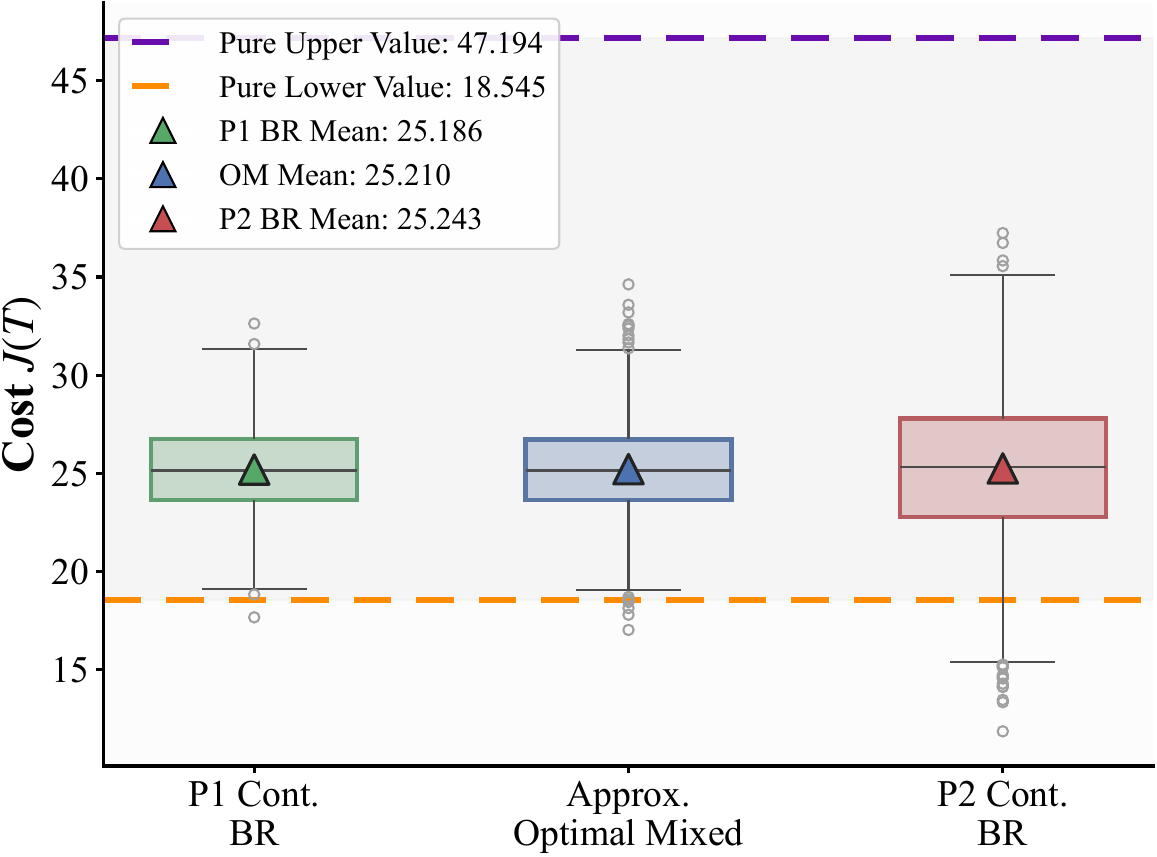}}
        \label{fig:suboptimality}
    }
    \caption{Weak approximation and suboptimality gap under different designs for the diffusion term of $\mathbf{\tilde{G}}$.}
    \label{fig:performance_certification}
\end{figure*}

To quantify the state distributional approximation, we measure the 1-Wasserstein distance, $W_1(\mu_{DG}(t), \mu_{SDG}(t))$, between the empirical state distributions of the true game and the auxiliary SDG at each time step. As shown in Fig. \ref{fig:wasserstein_compare}, the delay-scaled diffusion design ($\sigma_3$) successfully bounds the Wasserstein distance globally, demonstrating superior weak approximation capability. This structural alignment is further validated in Fig. \ref{fig:cost_boxplot}, where the expectation of the terminal cost $J(T)$ of the original DG almost perfectly coincides with the designed SDG under $\sigma_3$.

\subsection{Suboptimality and Pure Strategy Value Gap}
To demonstrate the advantage of the near-optimal mixed strategy and certify its suboptimality gap, we consider three strategy configurations: (1) Near Optimal Mixed (OM), (2) Player 1 Best Response (P1 BR) against Player 2's OM, and (3) Player 2 Best Response (P2 BR) against Player 1's OM.

Furthermore, we explicitly compute the theoretical \textit{Pure Upper Value} (P1 exploited unconditionally) and \textit{Pure Lower Value} (P2 exploited unconditionally) by solving the first-order HJI equation of the game without commitment delays. These bounds are depicted as the shaded gap in Fig. \ref{fig:suboptimality}.

The results show that the Monte Carlo mean of the proposed OM method lies strictly within the pure strategy gap. It empirically verifies that relying solely on deterministic pure strategies leads to deterministic exploitation, whereas our mixed strategy guarantees a robust expected performance strictly bounded away from the worst-case pure strategy outcomes.

\begin{figure*}[h]
    \centering
    \subfigure[weak approx: state Wasserstein distance]{
        \includegraphics[width = 0.31\textwidth]{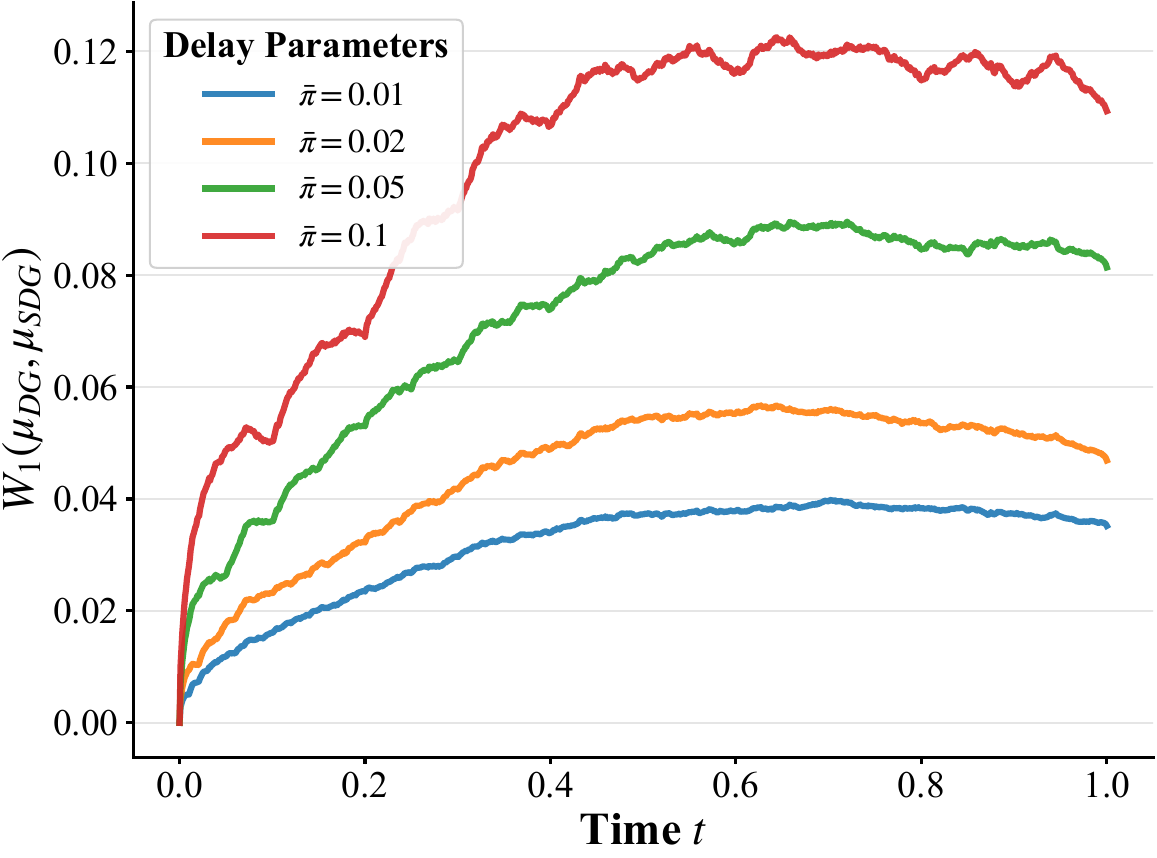}
        \label{fig:sens_wasserstein}
    }
    \subfigure[weak approx: cost expectation alignment]{
        \includegraphics[width = 0.31\textwidth]{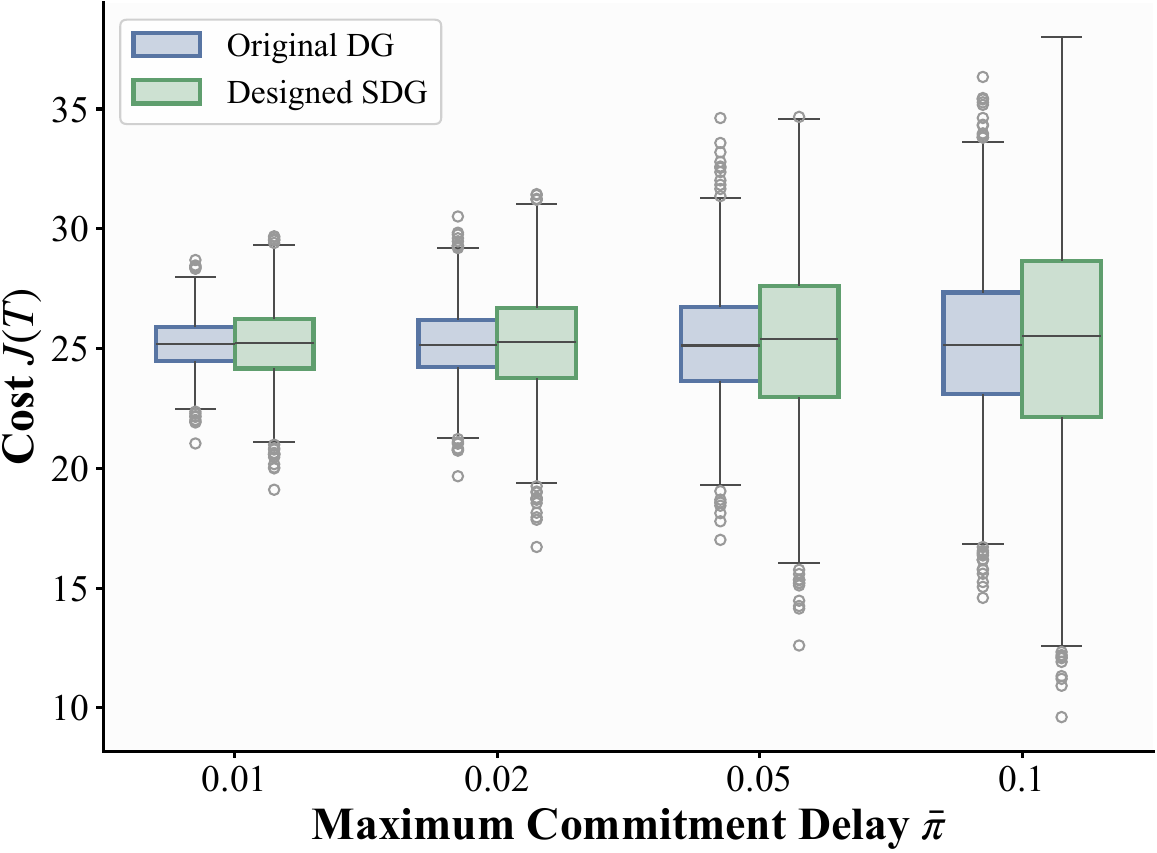}
        \label{fig:sens_cost}
    }
    \subfigure[suboptimality \& pure strategy value gap]{
        \includegraphics[width = 0.31\textwidth]{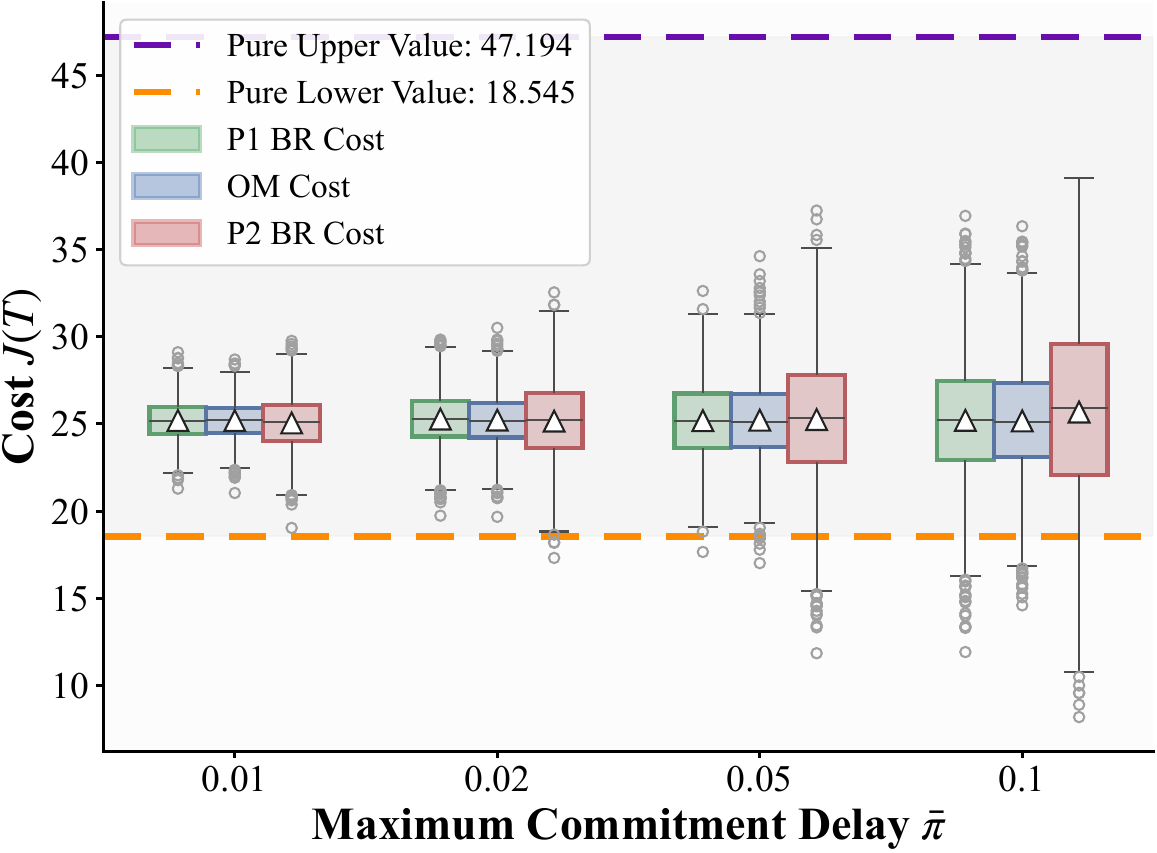}
        \label{fig:sens_suboptimality}
    }
    \caption{Weak approximation and suboptimality gap under different $\bar\pi$.}
    \label{fig:sens_performance_certification}
\end{figure*}
\subsection{Sensitivity to Commitment Delay $\pi$}
Finally, to provide a definitive verification of Theorem \ref{thm:general_performance}, we conduct a sensitivity analysis on the commitment delay by varying $\bar{\pi} \in \{0.10, 0.05, 0.02, 0.01\}$ under the diffusion design $\sigma_3$.

As shown in Fig. \ref{fig:sens_wasserstein} and Fig. \ref{fig:sens_cost}, the global Wasserstein distance between the state distributions and the expectation difference bewteen the cost decrease monotonically as $\bar{\pi}$ shrinks, validating the $\mathcal{O}(\bar{\pi})$ distributional convergence. More importantly, the grouped boxplots in Fig. \ref{fig:sens_suboptimality} explicitly reveal a ``funnel effect'': as the delay $\bar{\pi} \to 0$, the suboptimality gaps tightly shrink, and the expected costs of the Best Responses and the Optimal Mixed strategy strictly converge toward the center of the pure strategy value bounds. This observation powerfully substantiates our theoretical claim that the proposed discretization framework universally bounds both the value approximation error and the suboptimality gap by $\mathcal{O}(\bar{\pi})$.
\begin{remark}[Justification on surrogate best response and $\mathcal{O}(\bar{\pi})$ inversion]
\label{rem:surrogate_br}
In the suboptimality verification (Fig. \ref{fig:sens_suboptimality}), one may observe that the Monte Carlo evaluation of Player 1's best response cost, $J(\alpha_{BR}, \beta^d)$, is not strictly less than the approximated optimal mixed strategy cost, $J(\alpha^d, \beta^d)$. Occasionally, $J(\alpha_{BR}, \beta^d)$ might even be slightly greater. 

This phenomenon is theoretically sound and stems from the intractability of the exact best response. Computing the strict infimum $\inf_{\alpha \in \mathcal{A}_m^{\pi}} J(\alpha, \beta^d)$ over the original game $\mathbf{G}$ is computationally prohibitive, as it requires traversing an exponentially growing scenario tree generated by the opponent's delayed randomizations. Consequently, we compute a \textit{surrogate best response} by optimizing over the continuous-time auxiliary SDG $\mathbf{\tilde{G}}$, i.e., solving for the optimal continuous control against the opponent's approximated dynamics $\tilde{J}(\cdot, \phi_2(\beta^d))$. 

Let $\tilde{u}_{BR}$ be the optimal surrogate response. By definition, it guarantees $\tilde{J}(\tilde{u}_{BR}, \tilde{v}^*) \leq \tilde{J}(\tilde{u}^*, \tilde{v}^*)$ in the auxiliary SDG. However, mapping these strategies back to the true game $\mathbf{G}$ incurs weak approximation errors due to the ZOH execution. Using Theorem \ref{thm:value_approximation}, we have:
\begin{align*}
    &J(\alpha_{BR}, \beta^d) - J(\alpha^d, \beta^d)\\ 
    &\leq \left( \tilde{J}(\tilde{u}_{BR}, \tilde{v}^*) + \mathcal{O}(\bar{\pi}) \right) - \left( \tilde{J}(\tilde{u}^*, \tilde{v}^*) - \mathcal{O}(\bar{\pi}) \right)\leq \mathcal{O}(\bar{\pi}).
\end{align*}
Therefore, the slight inversion where $J(\alpha_{BR}, \beta^d) \gtrsim J(\alpha^d, \beta^d)$ is strictly bounded by the $\mathcal{O}(\bar{\pi})$ gap. Rather than refuting the methodology, this bounded discrepancy practically validates the tightness of our weak approximation design under the given commitment delay $\bar{\pi}$.
\end{remark}

\section{Conclusion}
This paper has investigated the challenging problem of computing optimal mixed strategies for zero-sum differential games. We introduced a weak approximation framework that establishes a stochastic equivalent representation of the original game, which ensures the existence of a game value and provides a certified approximation of the state dynamics and cost expectations. To facilitate implementation, we further developed a control-space discretization algorithm, transforming the abstract measure-theoretic optimization into solvable local matrix games. 
Theoretical analysis demonstrates that the weak approximation error is strictly of order $\mathcal{O}(\bar\pi)$ with respect to the commitment delay, while explicit analytical upper bounds are provided to certify the strategy suboptimality gaps. Numerical simulations further validate the theoretical convergence bounds and demonstrate the advantage of the synthesized mixed strategies in preventing deterministic exploitation. 

Extending this methodology to $N$-player general-sum differential games is a promising direction. Our core framework of mapping continuous-time mixed strategies to step-wise pure strategies remains valid. However, the equilibrium analysis must generalize from a unique saddle point to general Nash Equilibria. Computationally, this requires transitioning from a single Hamilton-Jacobi-Isaacs equation to a system of coupled HJI equations. Fortunately, our discrete probability parameterization naturally facilitates this extension by leveraging algorithms developed for standard extensive-form games. We hope this framework motivates further investigations in these broader settings.
\bibliographystyle{IEEEtran}
\bibliography{refs.bib}

@article{flemingExistenceValueFunctions1989b,
  title = {On the Existence of Value Functions of Two-Player, Zero-Sum Stochastic Differential Games},
  author = {Fleming, Wendell H. and Souganidis, Panagiotis E.},
  year = 1989,
  journal = {Indiana University Mathematics Journal},
  volume = {38},
  number = {2},
  pages = {293--314},
  publisher = {JSTOR}
}

@incollection{aumann28MixedBehavior1964,
  title = {28. Mixed and Behavior Strategies in Infinite Extensive Games},
  booktitle = {Advances in Game Theory. (AM-52)},
  author = {Aumann, Robert J.},
  year = {1964},
  month = dec,
  pages = {627--650},
  publisher = {Princeton University Press},
  urldate = {2024-04-09},
  isbn = {978-1-4008-8201-4}
}

@article{aumannBorelStructuresFunction1961,
  title = {Borel Structures for Function Spaces},
  author = {Aumann, Robert J.},
  year = {1961},
  journal = {Illinois Journal of Mathematics},
  volume = {5},
  number = {4},
  pages = {614--630},
  publisher = {Duke University Press},
  urldate = {2024-04-09}
}

@book{basarDynamicNoncooperativeGame1998,
  title = {Dynamic Noncooperative Game Theory},
  author = {Ba{\c s}ar, Tamer and Olsder, Geert Jan},
  year = {1998},
  publisher = {SIAM}
}

@book{basarHandbookDynamicGame2018,
  title = {Handbook of Dynamic Game Theory},
  author = {Ba{\c s}ar, Tamer and Zaccour, Georges},
  year = {2018},
  month = jun,
  publisher = {Springer International Publishing},
  googlebooks = {WOCtDAEACAAJ},
  isbn = {978-3-319-44373-7},
  langid = {english}
}

@incollection{berkovitzRelaxedControls2012,
  title = {Relaxed Controls},
  booktitle = {Nonlinear Optimal Control Theory},
  author = {Berkovitz, Leonard David and Medhin, Negash G.},
  year = {2012},
  publisher = {{Chapman and Hall/CRC}},
  isbn = {978-0-429-09810-9}
}

@article{buckdahnValueFunctionDifferential2013,
  title = {Value Function of Differential Games without Isaacs Conditions. An Approach with Nonanticipative Mixed Strategies},
  author = {Buckdahn, Rainer and Li, Juan and Quincampoix, Marc},
  year = {2013},
  month = nov,
  journal = {International Journal of Game Theory},
  volume = {42},
  number = {4},
  pages = {989--1020},
  issn = {1432-1270},
  urldate = {2023-03-04},
  langid = {english}
}

@article{buckdahnValueMixedStrategies2014,
  title = {Value in Mixed Strategies for Zero-Sum Stochastic Differential Games Without Isaacs Condition},
  author = {Buckdahn, Rainer and Li, Juan and Quincampoix, Marc},
  year = {2014},
  journal = {The Annals of Probability},
  volume = {42},
  number = {4},
  eprint = {42920505},
  eprinttype = {jstor},
  pages = {1724--1768},
  publisher = {Institute of Mathematical Statistics},
  issn = {0091-1798},
  urldate = {2023-03-29}
}

@article{cardaliaguetDifferentialGamesAsymmetric2007,
  title = {Differential Games with Asymmetric Information},
  author = {Cardaliaguet, P.},
  year = {2007},
  month = jan,
  journal = {SIAM Journal on Control and Optimization},
  volume = {46},
  number = {3},
  pages = {816--838},
  publisher = {{Society for Industrial and Applied Mathematics}},
  issn = {0363-0129},
  urldate = {2023-02-09}
}

@article{cardaliaguetPureRandomStrategies2014,
  title = {Pure and Random Strategies in Differential Game with Incomplete Informations},
  author = {Cardaliaguet, Pierre and Jimenez, Chlo{\'e} and Quincampoix, Marc},
  year = {2014},
  journal = {Journal of Dynamics and Games},
  volume = {1},
  number = {3},
  pages = {363--375},
  publisher = {{Journal of Dynamics and Games}},
  issn = {2164-6066},
  urldate = {2023-02-09},
  copyright = {http://creativecommons.org/licenses/by/3.0/},
  langid = {english}
}

@article{coxControlledMeasurevaluedMartingales2024,
  title = {Controlled Measure-Valued Martingales: A Viscosity Solution Approach},
  shorttitle = {Controlled Measure-Valued Martingales},
  author = {Cox, Alexander M. G. and K{\"a}llblad, Sigrid and Larsson, Martin and {Svaluto-Ferro}, Sara},
  year = {2024},
  month = apr,
  journal = {The Annals of Applied Probability},
  volume = {34},
  number = {2},
  pages = {1987--2035},
  publisher = {Institute of Mathematical Statistics},
  issn = {1050-5164, 2168-8737},
  urldate = {2024-05-04}
}

@article{elliottExistenceValueStochastic1976,
  title = {The Existence of Value in Stochastic Differential Games},
  author = {Elliott, Robert},
  year = {1976},
  month = jan,
  journal = {SIAM Journal on Control and Optimization},
  volume = {14},
  number = {1},
  pages = {85--94},
  publisher = {{Society for Industrial and Applied Mathematics}},
  issn = {0363-0129},
  urldate = {2023-03-28}
}

@book{engwerdaLQDynamicOptimization2005,
  title = {LQ Dynamic Optimization and Differential Games},
  author = {Engwerda, Jacob},
  year = {2005},
  month = jun,
  edition = {1st edition},
  publisher = {Wiley},
  address = {Chicester, West Sussex, England Hoboken, NJ},
  isbn = {978-0-470-01524-7},
  langid = {english}
}

@article{evansDifferentialGamesRepresentation1984c,
  title = {Differential Games and Representation Formulas for Solutions of Hamilton-Jacobi-Isaacs Equations},
  author = {Evans, Lawrence C. and Souganidis, Panagiotis E.},
  year = {1984},
  journal = {Indiana University mathematics journal},
  volume = {33},
  number = {5},
  eprint = {45010271},
  eprinttype = {jstor},
  pages = {773--797},
  publisher = {JSTOR},
  urldate = {2025-03-25}
}

@inproceedings{evenContinuizedAccelerationsDeterministic2021,
  title = {Continuized Accelerations of Deterministic and Stochastic Gradient Descents, and of Gossip Algorithms},
  booktitle = {Advances in Neural Information Processing Systems},
  author = {Even, Mathieu and Berthier, Rapha{\"e}l and Bach, Francis and Flammarion, Nicolas and Hendrikx, Hadrien and Gaillard, Pierre and Massouli{\'e}, Laurent and Taylor, Adrien},
  year = {2021},
  volume = {34},
  pages = {28054--28066},
  urldate = {2025-03-25}
}

@book{flemingControlledMarkovProcesses2006,
  title = {Controlled Markov Processes and Viscosity Solutions},
  author = {Fleming, Wendell H. and Soner, Halil Mete},
  year = {2006},
  month = feb,
  publisher = {Springer Science \& Business Media},
  googlebooks = {4Bjz2iWmLyQC},
  isbn = {978-0-387-31071-8},
  langid = {english}
}

@article{flemingMixedStrategiesDeterministic2017a,
  title = {Mixed Strategies for Deterministic Differential Games},
  author = {Fleming, Wendell H. and {Hernandez-Hernandez}, Daniel},
  year = {2017},
  month = jun,
  journal = {Communications on Stochastic Analysis},
  volume = {11},
  number = {2},
  issn = {0973-9599},
  urldate = {2023-05-17},
  langid = {english}
}

@inproceedings{goodfellowGenerativeAdversarialNets2014,
  title = {Generative Adversarial Nets},
  booktitle = {Advances in Neural Information Processing Systems},
  author = {Goodfellow, Ian and {Pouget-Abadie}, Jean and Mirza, Mehdi and Xu, Bing and {Warde-Farley}, David and Ozair, Sherjil and Courville, Aaron and Bengio, Yoshua},
  year = {2014},
  volume = {27},
  urldate = {2024-05-17}
}

@article{harrisExistenceSubgameperfectEquilibrium1995,
  title = {The Existence of Subgame-Perfect Equilibrium in Continuous Games with Almost Perfect Information: A Case for Public Randomization},
  shorttitle = {The Existence of Subgame-Perfect Equilibrium in Continuous Games with Almost Perfect Information},
  author = {Harris, Christopher and Reny, Philip and Robson, Arthur},
  year = {1995},
  journal = {Econometrica: Journal of the Econometric Society},
  eprint = {2171906},
  eprinttype = {jstor},
  pages = {507--544},
  publisher = {JSTOR},
  urldate = {2025-03-25}
}

@inproceedings{hespanhaProbabilisticPursuitevasionGames2000,
  title = {Probabilistic Pursuit-Evasion Games: A One-Step Nash Approach},
  shorttitle = {Probabilistic Pursuit-Evasion Games},
  booktitle = {Proceedings of the 39th IEEE Conference on Decision and Control},
  author = {Hespanha, J.P. and Prandini, M. and Sastry, S.},
  year = {2000},
  month = dec,
  volume = {3},
  pages = {2272-2277 vol.3},
  issn = {0191-2216}
}

@book{isaacsDifferentialGamesMathematical1999,
  title = {Differential Games: A Mathematical Theory with Applications to Warfare and Pursuit, Control and Optimization},
  shorttitle = {Differential Games},
  author = {Isaacs, Rufus},
  year = {1999},
  publisher = {Courier Corporation}
}

@article{islerRandomizedPursuitEvasionLocal2006,
  title = {Randomized Pursuit-Evasion with Local Visibility},
  author = {Isler, Volkan and Kannan, Sampath and Khanna, Sanjeev},
  year = {2006},
  month = jan,
  journal = {SIAM Journal on Discrete Mathematics},
  volume = {20},
  number = {1},
  pages = {26--41},
  publisher = {{Society for Industrial and Applied Mathematics}},
  issn = {0895-4801},
  urldate = {2023-06-19}
}

@article{islerRandomizedPursuitevasionPolygonal2005,
  title = {Randomized Pursuit-Evasion in a Polygonal Environment},
  author = {Isler, V. and Kannan, S. and Khanna, S.},
  year = {2005},
  month = oct,
  journal = {IEEE Transactions on Robotics},
  volume = {21},
  number = {5},
  pages = {875--884},
  issn = {1941-0468}
}

@article{kumarOptimalMixedStrategies1980,
  title = {Optimal Mixed Strategies in a Dynamic Game},
  author = {Kumar, P.},
  year = {1980},
  month = aug,
  journal = {IEEE Transactions on Automatic Control},
  volume = {25},
  number = {4},
  pages = {743--749},
  issn = {1558-2523},
  urldate = {2024-04-23}
}

@article{kunischOptimalControlUndamped2016,
  title = {Optimal Control of the Undamped Linear Wave Equation with Measure Valued Controls},
  author = {Kunisch, Karl and Trautmann, Philip and Vexler, Boris},
  year = {2016},
  month = jan,
  journal = {SIAM Journal on Control and Optimization},
  volume = {54},
  number = {3},
  pages = {1212--1244},
  publisher = {{Society for Industrial and Applied Mathematics}},
  issn = {0363-0129},
  urldate = {2024-05-04}
}

@article{latzAnalysisStochasticGradient2021a,
  title = {Analysis of Stochastic Gradient Descent in Continuous Time},
  author = {Latz, Jonas},
  year = {2021},
  month = may,
  journal = {Statistics and Computing},
  volume = {31},
  number = {4},
  pages = {39},
  issn = {1573-1375},
  urldate = {2024-05-03},
  langid = {english}
}

@article{levantChatteringAnalysis2010,
  title = {Chattering Analysis},
  author = {Levant, Arie},
  year = {2010},
  month = jun,
  journal = {IEEE Transactions on Automatic Control},
  volume = {55},
  number = {6},
  pages = {1380--1389},
  issn = {1558-2523},
  urldate = {2025-03-19}
}

@inproceedings{liStochasticModifiedEquations2017,
  title = {Stochastic Modified Equations and Adaptive Stochastic Gradient Algorithms},
  booktitle = {International Conference on Machine Learning},
  author = {Li, Qianxiao and Tai, Cheng and Weinan, E.},
  year = {2017},
  pages = {2101--2110},
  publisher = {PMLR}
}

@article{liStochasticModifiedEquations2019,
  title = {Stochastic Modified Equations and Dynamics of Stochastic Gradient Algorithms i: Mathematical Foundations},
  shorttitle = {Stochastic Modified Equations and Dynamics of Stochastic Gradient Algorithms i},
  author = {Li, Qianxiao and Tai, Cheng and Weinan, E.},
  year = {2019},
  journal = {The Journal of Machine Learning Research},
  volume = {20},
  number = {1},
  pages = {1474--1520},
  publisher = {JMLR. org}
}

@inproceedings{martinFindingMixedstrategyEquilibria2023,
  title = {Finding Mixed-Strategy Equilibria of Continuous-Action Games without Gradients Using Randomized Policy Networks},
  booktitle = {Proceedings of the Thirty-Second International Joint Conference on Artificial Intelligence},
  author = {Martin, Carlos and Sandholm, Tuomas},
  year = {2023},
  month = aug,
  series = {IJCAI '23},
  pages = {2844--2852},
  urldate = {2024-05-10},
  isbn = {978-1-956792-03-4}
}

@article{milshteinWeakApproximationSolutions1986a,
  title = {Weak Approximation of Solutions of Systems of Stochastic Differential Equations},
  author = {Mil'shtein, G. N.},
  year = {1986},
  journal = {Theory of Probability \& Its Applications},
  volume = {30},
  number = {4},
  pages = {750--766},
  publisher = {SIAM}
}

@incollection{nashjrNoncooperativeGames1996,
  title = {Non-Cooperative Games},
  booktitle = {Essays on Game Theory},
  author = {Nash Jr, John},
  year = {1996},
  pages = {22--33},
  publisher = {Edward Elgar Publishing}
}

@article{nikaidoNeumannsMinimaxTheorem1954,
  title = {On von Neumann's Minimax Theorem},
  author = {Nikaid{\^o}, Hukukane},
  year = {1954},
  journal = {Pacific J. Math},
  volume = {4},
  number = {1},
  pages = {65--72},
  urldate = {2024-04-25}
}

@book{osborneCourseGameTheory1994,
  title = {A Course in Game Theory},
  author = {Osborne, Martin J. and Rubinstein, Ariel},
  year = {1994},
  month = jul,
  publisher = {The MIT Press},
  address = {Cambridge, Mass},
  isbn = {978-0-262-65040-3},
  langid = {english}
}

@article{perkinsMixedStrategyLearningContinuous2017a,
  title = {Mixed-Strategy Learning With Continuous Action Sets},
  author = {Perkins, Steven and Mertikopoulos, Panayotis and Leslie, David S.},
  year = {2017},
  month = jan,
  journal = {IEEE Transactions on Automatic Control},
  volume = {62},
  number = {1},
  pages = {379--384},
  issn = {1558-2523},
  urldate = {2024-05-06}
}

@inproceedings{petersLearningMixedStrategies2022,
  title = {Learning Mixed Strategies in Trajectory Games},
  booktitle = {Robotics: Science and Systems XVIII},
  author = {Peters, Lasse and {Fridovich-Keil}, David and Ferranti, Laura and Stachniss, Cyrill and {Alonso-Mora}, Javier and Laine, Forrest},
  year = {2022},
  month = jun,
  publisher = {{Robotics: Science and Systems Foundation}},
  urldate = {2024-04-23},
  isbn = {978-0-9923747-8-5},
  langid = {english}
}

@article{roxinAxiomaticApproachDifferential1969,
  title = {Axiomatic Approach in Differential Games},
  author = {Roxin, Emilio},
  year = {1969},
  month = mar,
  journal = {Journal of Optimization Theory and Applications},
  volume = {3},
  number = {3},
  pages = {153--163},
  issn = {1573-2878},
  urldate = {2025-03-25},
  langid = {english}
}

@article{sirbuStochasticPerronsMethod2014,
  title = {Stochastic Perron's Method and Elementary Strategies for Zero-Sum Differential Games},
  author = {S{\^i}rbu, Mihai},
  year = {2014},
  month = jan,
  journal = {SIAM Journal on Control and Optimization},
  volume = {52},
  number = {3},
  pages = {1693--1711},
  publisher = {{Society for Industrial and Applied Mathematics}},
  issn = {0363-0129},
  urldate = {2024-04-07}
}

@article{varaiyaExistenceSaddlePoints1969,
  title = {Existence of Saddle Points in Differential Games},
  author = {Varaiya, P. and Lin, Jiguan},
  year = {1969},
  month = feb,
  journal = {SIAM Journal on Control},
  volume = {7},
  number = {1},
  pages = {141--157},
  issn = {0036-1402},
  urldate = {2025-03-25},
  langid = {english}
}

@article{vidalProbabilisticPursuitevasionGames2002a,
  title = {Probabilistic Pursuit-Evasion Games: Theory, Implementation, and Experimental Evaluation},
  shorttitle = {Probabilistic Pursuit-Evasion Games},
  author = {Vidal, R. and Shakernia, O. and Kim, H.J. and Shim, D.H. and Sastry, S.},
  year = {2002},
  month = oct,
  journal = {IEEE Transactions on Robotics and Automation},
  volume = {18},
  number = {5},
  pages = {662--669},
  issn = {2374-958X}
}

@article{wangPontryaginsMaximumPrinciple2010,
  title = {A Pontryagin's Maximum Principle for Non-Zero Sum Differential Games of BSDEs with Applications},
  author = {Wang, Guangchen and Yu, Zhiyong},
  year = {2010},
  month = jul,
  journal = {IEEE Transactions on Automatic Control},
  volume = {55},
  number = {7},
  pages = {1742--1747},
  issn = {1558-2523}
}

@inproceedings{xuDifferentialGameMixed2023b,
  title = {Differential Game with Mixed Strategies: A Weak Approximation Approach},
  shorttitle = {Differential Game with Mixed Strategies},
  booktitle = {2023 62nd IEEE Conference on Decision and Control (CDC)},
  author = {Xu, Tao and Xi, Wang and He, Jianping},
  year = {2023},
  month = dec,
  pages = {5216--5221},
  issn = {2576-2370},
  urldate = {2025-03-25}
}

@article{dou2019finding,
  title={Finding mixed strategy nash equilibrium for continuous games through deep learning},
  author={Dou, Zehao and Yan, Xiang and Wang, Dongge and Deng, Xiaotie},
  journal={arXiv preprint arXiv:1910.12075},
  year={2019}
}

\begin{appendices}
\section{Proof of Theorem \ref{thm:SE}} \label{app:thm:SE}
\begin{proof}
Our proof proceeds based on the principle of dynamic programming. Notice that the game value functions, presented in Definition \ref{def:value-function}, are defined on the initial state and time without concerning the information filtration generated by the mixed strategies. Nevertheless, to further define the value function at any time $t\in(t_0,T]$ during the game, one should notice that the value functions is measurable with respect to $\mathcal{F}_t$. 
Therefore, the first step of our proof is to show that the value functions are not influenced by $\mathcal{F}_t$. First, for $k\in\{0,\ldots,N-1\}$, we use the essential supremum and essential infimum over $\mathcal{F}_{k}$ to define value functions at $t_k$, i.e.,
\begin{equation*}
\begin{aligned}
    &V_{m,-}^\pi(t_{k},x(t_k)) \\
    & \qquad\qquad \triangleq\operatorname{esssup}\limits_{\beta\in\mathcal{B}_m^\pi}\operatorname{essinf}\limits_{\alpha\in\mathcal{A}_m^\pi} J(t_{k},x(t_k),\alpha,\beta)\!\mid\!\mathcal{F}_{k},\\
    &V_{m,+}^\pi(t_{k},x(t_k)) \\
    & \qquad\qquad\triangleq \operatorname{essinf}\limits_{\alpha\in\mathcal{A}_m^\pi}\operatorname{esssup}\limits_{\beta\in\mathcal{B}_m^\pi} J(t_{k},x(t_k),\alpha,\beta)\!\mid\!\mathcal{F}_{k}.
\end{aligned}
\end{equation*}
It follows that
\begin{equation*}
    \begin{aligned}
                           & V_{m,-}^\pi(t_{N-1},x(t_{N-1}))                                                                                                                                                                \\
        \overset{(i)}{=}   & \operatorname{esssup}_{\beta\in \mathcal{B}_m^{\pi}}\operatorname{essinf}_{\alpha\in\mathcal{A}_m^{\pi}} \mathbb{E}_{\omega}\left\{ \int_{t_{N-1}}^{t_{N}}\!\!h(t,x(t),u^\omega,v^\omega)dt \right.                                                          \\
                           & \qquad \qquad \qquad \qquad+ V_{m,-}^\pi(t_{N}, x(t_{N})) \mid \mathcal{F}_{N-1} \big\}                                                                                                                               \\
        \overset{(ii)}{=}  & \sup_{\nu\in \mathcal{P}_{\mathcal{C}}(V)}\inf_{\mu\in \mathcal{P}_{\mathcal{C}}(U)} \int_{V\times U}\left\{\int_{t_{N-1}}^{t_{N}}\!\!h(t,x(t),u,v)ds \right. \\
                           & \qquad \qquad \qquad \qquad \qquad+V_{m,-}^\pi(t_{N}, x(t_{N})) \big\}\; \!\!\nu(dv)\mu(du). 
    \end{aligned}
\end{equation*}
Subject to the convex constraint set $\mathcal{C}$, $\mathcal{P}_{\mathcal{C}}(U)$ and $\mathcal{P}_{\mathcal{C}}(V)$ are convex 
subsets of the totalities of all regular Borel probability containing probability measures on $U$ and $V$. For equation $(i)$, $(u^\omega, v^\omega)$ is a unique pair of control functions belonging to $\operatorname{ST}(t_0,T;\pi, U) \times \operatorname{ST}(t_0,T;\pi,V)$ such that $\alpha(w,v^\omega) \equiv u^\omega$ and $\beta(w,u^\omega) \equiv v^\omega$, for any $\omega\in\Omega$. 
For equation $(ii)$, on the one hand, notice that $V_{m,-}^\pi(t_{N},x(t_{N})) = g(x(t_{N}))$, which only depends on $x(t_{N-1})$ and $\mathcal{G}_{N}$, and is independent of $\mathcal{F}_{N-1}$.
On the other hand, since the mixed strategies determine the probability distributions of $u^\omega_{\mid_{[t_{N-1},t_N]}}$ and $v^\omega_{\mid_{[t_{N-1},t_N]}}$, the essential supremum over $\mathcal{B}_m^\pi$ and essential infimum over $\mathcal{A}_m^\pi$ are equivalent to the supremum over $\mathcal{P}_{\mathcal{C}}(V)$ and infimum over $\mathcal{P}_{\mathcal{C}}(U)$, resp. Moreover, the measure spaces $\mathcal{P}_{\mathcal{C}}(V)$ and $\mathcal{P}_{\mathcal{C}}(U)$ are independent of $\mathcal{F}_{N-1}$.
Therefore, one gets that $V_{m,-}^\pi(t_{N-1}, x(t_{N-1}))$ does not depend on $\mathcal{F}_{N-1}$. 
By mathematical induction, it follows that $V_{m,-}^\pi(t_{k-1},x)$ does not depend on $\mathcal{F}_{k-1}$ for $k\in\{1,\ldots,N\}$. Similarly, this conclusion also applies to the upper value functions.

The second step of the proof is to show that the upper and lower values coincide, which is then an obvious result because
	\begin{equation*}
		\begin{aligned}
			                 & V_{m,-}^\pi(t_{k-1},x(t_{k-1}))                                                                                                                                                                            \\
			=                & \sup_{\nu\in \mathcal{P}_{\mathcal{C}}(V)}\inf_{\mu\in \mathcal{P}_{\mathcal{C}}(U)}  \int_{V\times U}\left\{\int_{t_{k-1}}^{t_{k}}\!\!h(t,x(t),u,v)ds\right. \\
			                 & \qquad \qquad \qquad \qquad \qquad+\left.V_{m,-}^\pi(t_{k}, x(t_{k}))\right\}\; \nu(dv)\mu(du)                                                                                   \\
			\overset{(i)}{=} & \inf_{\mu\in \mathcal{P}_{\mathcal{C}}(U)} \sup_{\nu\in \mathcal{P}_{\mathcal{C}}(V)} \int_{U\times V}\left\{\int_{t_{k-1}}^{t_{k}}h(t,x(t),u,v)dt\right.     \\
			                 & \left.\qquad \qquad \qquad \qquad \qquad+V_{m,-}^\pi(t_{k}, x(t_{k}))\right\}\; \mu(du)\nu(dv)                                                                                   \\
			=                & \;V_{m,+}^\pi(t_{k-1},x(t_{k-1})),
		\end{aligned}
	\end{equation*}
	where $(i)$ holds because the general Von Neumann minimax theorem guarantees the existence of a saddle-point \cite{nikaidoNeumannsMinimaxTheorem1954}.
\end{proof}

\section{Proof of Theorem \ref{thm:value_approximation}} \label{app:thm:value_approximation}
\begin{proof}

    For ease of notation, we omit the initial conditions $(t_0,x_0)$ in $J$ and $\tilde{J}$ throughout this proof without introducing ambiguity. \\
    \underline{\textit{Step 1}.} Define $a = J(\alpha^*, \beta^d)$, $b = J(\alpha^d, \beta^*)$, $c = \tilde{J}(\phi_1(\alpha^*), \tilde{\beta}^*)$, and $d = \tilde{J}(\tilde{\alpha}^*, \phi_2(\beta^*))$. By the $n$-th order weak approximation between $\mathbf{G}$ and $\mathbf{\tilde{G}}$, we obtain $|a - \tilde{J}(\phi_1(\alpha^*), \phi_2(\beta^d))| = \mathcal{O}(\bar{\pi}^n)$ and $|b - \tilde{J}(\phi_1(\alpha^d), \phi_2(\beta^*))| = \mathcal{O}(\bar{\pi}^n)$. Additionally, by definition, $|d - \tilde{J}(\phi_1(\alpha^d), \phi_2(\beta^*))| = \epsilon_1$ and $|c - \tilde{J}(\phi_1(\alpha^*), \phi_2(\beta^d))| = \epsilon_2$. Consequently, we can establish the following relations:
    \begin{equation*}
        \begin{array}{ccccc}
        a = J(\alpha^*, \beta^d) & \le & J(\alpha^*, \beta^*) & \le & J(\alpha^d, \beta^*) = b \\
        \Big\updownarrow {\scriptstyle \mathcal{O}(\bar{\pi}^n)} & & & & \Big\updownarrow {\scriptstyle \mathcal{O}(\bar{\pi}^n)} \\
        \tilde{J}(\phi_1(\alpha^*), \phi_2(\beta^d)) & & & & \tilde{J}(\phi_1(\alpha^d), \phi_2(\beta^*)) \\
        \Big\updownarrow {\scriptstyle \epsilon_2} & & & & \Big\updownarrow {\scriptstyle \epsilon_1} \\
        c = \tilde{J}(\phi_1(\alpha^*), \tilde{\beta}^*) & \ge & \tilde{J}(\tilde{\alpha}^*, \tilde{\beta}^*) & \ge & \tilde{J}(\tilde{\alpha}^*, \phi_2(\beta^*)) = d.
        \end{array}
    \end{equation*}
    
    \underline{\textit{Step 2}.} Next, we prove that 
    \begin{equation}\label{eq:inequality_lem}
        \max(|a-d|, |b-c|) \leq \max(|a-c|, |b-d|).
    \end{equation}
    Recall the algebraic identity for the maximum of two absolute values:
    \begin{equation}
        \max(|X|, |Y|) = \frac{|X + Y| + |X - Y|}{2}, \quad \forall X, Y \in \mathbb{R}. \label{eq:max_abs_identity}
    \end{equation}
    Applying \eqref{eq:max_abs_identity} to the left-hand side (LHS) of the proposed inequality, we obtain:
    \begin{align}
        \text{LHS} &= \frac{|(a-d) + (b-c)| + |(a-d) - (b-c)|}{2} \nonumber \\
        &= \frac{1}{2} \Big( |(a+b) - (c+d)| + |(a-b) + (c-d)| \Big). \label{eq:lhs_bound}
    \end{align}
    Similarly, applying \eqref{eq:max_abs_identity} to the right-hand side (RHS) yields:
    \begin{align}
        \text{RHS} &= \frac{|(a-c) + (b-d)| + |(a-c) - (b-d)|}{2} \nonumber \\
        &= \frac{1}{2} \Big( |(a+b) - (c+d)| + |(a-b) - (c-d)| \Big). \label{eq:rhs_bound}
    \end{align}
    Comparing \eqref{eq:lhs_bound} and \eqref{eq:rhs_bound}, we observe that the first term, $|(a+b) - (c+d)|$, is identical on both sides. Therefore, to establish that $\text{LHS} \leq \text{RHS}$, it is mathematically sufficient to prove that the second term of the LHS is bounded by the second term of the RHS:
    \begin{equation}
        |(a-b) + (c-d)| \leq |(a-b) - (c-d)|. \label{eq:reduced_abs_ineq}
    \end{equation}

    To evaluate \eqref{eq:reduced_abs_ineq}, we introduce the variables $x \triangleq a - b$ and $y \triangleq c - d$. By hypothesis, $a \leq b$, which implies $x \leq 0$. Furthermore, $c \geq d$, which implies $y \geq 0$. Substituting these definitions into \eqref{eq:reduced_abs_ineq} gives:
    \begin{equation}
        |x + y| \leq |x - y|. \label{eq:xy_ineq}
    \end{equation}
    Given that $x \leq 0$ and $y \geq 0$, their product $xy$ is strictly non-positive ($xy \leq 0$), which implies $(x + y)^2 \leq (x - y)^2$. Thus, the condition in \eqref{eq:xy_ineq} is inherently satisfied. Since this inequality holds unconditionally under the given assumptions, the original statement in \eqref{eq:inequality_lem} is proven.

    \underline{\textit{Step 3}.} Finally, applying the triangle inequality and combining the established bounds, we conclude:
    \begin{equation*}
        \begin{aligned}
            &|V_m^\pi(t_0,x_0)-\tilde{V}_p(t_0,x_0)| = |J(\alpha^*, \beta^*) - \tilde{J}(\tilde{\alpha}^*, \tilde{\beta}^*)| \\
            \leq & \max(|a-d|, |b-c|) \leq \max(|a-c|, |b-d|) \\
            \leq & \max(\epsilon_1, \epsilon_2) + \mathcal{O}(\bar{\pi}^n).
        \end{aligned}
    \end{equation*}
\end{proof}

\section{Proof of Theorem \ref{thm:suboptimality_gap}}\label{app:thm:suboptimality_gap}
\begin{proof}
    For ease of notation, we omit $(t_0,x_0)$ in $J$ and $\tilde{J}$. We present the proof for the first inequality; the second follows symmetrically.

    \underline{\textit{Step 1}.} Consider strategies $\alpha^{(1)}, \alpha^{(2)} \in \mathcal{A}_m^\pi$ such that $\alpha^{(1)} \in \arg\min_{\alpha\in\mathcal{A}_m^{\pi}}J(\alpha, \beta^d)$ and $\phi_1(\alpha^{(2)}) = \tilde{\alpha}^{br}$. Let $a = J(\alpha^{(1)}, \beta^d)$, $b^\prime = J(\alpha^{(2)}, \beta^d)$, $c^\prime = \tilde{J}(\phi_1(\alpha^{(1)}), \phi_2(\beta^d))$, and $d^\prime = \tilde{J}(\tilde{\alpha}^{br}, \phi_2(\beta^d))$. According to the $n$-th order weak approximation between $\mathbf{G}$ and $\mathbf{\tilde{G}}$, we have $|a - c^\prime| = \mathcal{O}(\bar{\pi}^n)$ and $|b^\prime - d^\prime| = \mathcal{O}(\bar{\pi}^n)$. Consequently, the following relations hold:
    \begin{equation*}
        \begin{array}{ccc}
            a = J(\alpha^{(1)}, \beta^d) & \leq & J(\alpha^{(2)}, \beta^d) = b^\prime \\
            \Big\updownarrow {\scriptstyle \mathcal{O}(\bar{\pi}^n)} & & \Big\updownarrow {\scriptstyle \mathcal{O}(\bar{\pi}^n)}\\
            c^\prime = \tilde{J}(\phi_1(\alpha^{(1)}), \phi_2(\beta^d)) & \geq & \tilde{J}(\tilde{\alpha}^{br}, \phi_2(\beta^d)) = d^\prime.
        \end{array}
    \end{equation*}
    It follows that $|a-d^\prime| \leq \max(|a-d^\prime|, |b^\prime-c^\prime|) \leq \max(|a-c^\prime|, |b^\prime-d^\prime|) = \mathcal{O}(\bar{\pi}^n)$.

    \underline{\textit{Step 2}.} Let $b = J(\alpha^d, \beta^d)$, $c = \tilde{J}(\tilde{\alpha}^{br}, \beta^{*})$, and $d = \tilde{J}(\tilde{\alpha}^*, \tilde{\beta}^*)$. From Step 1, we have $|a - d^\prime| = \mathcal{O}(\bar{\pi}^n)$. Furthermore, $|b - \tilde{J}(\operatorname{ZOH}_{\pi}[\tilde{\alpha}^d], \operatorname{ZOH}_{\pi}[\tilde{\beta}^*])| = \mathcal{O}(\bar{\pi}^n)$ by the definition of the $n$-th order weak approximation. By definition of the error terms, $|d^\prime - c| = \epsilon_4$ and $|d - \tilde{J}(\operatorname{ZOH}_{\pi}[\tilde{\alpha}^*], \operatorname{ZOH}_{\pi}[\tilde{\beta}^*])| = \epsilon_3$. Hence, we establish:
    \begin{equation*}
        \begin{array}{ccc}
            a = J(\alpha^{(1)}, \beta^d) & \leq & J(\alpha^d, \beta^d) = b \\
            \Big\updownarrow {\scriptstyle \mathcal{O}(\bar{\pi}^n)} & & \Big\updownarrow {\scriptstyle \mathcal{O}(\bar{\pi}^n)}\\
            d^\prime = \tilde{J}(\tilde{\alpha}^{br}, \phi_2(\beta^d))& & \tilde{J}(\operatorname{ZOH}_{\pi}[\tilde{\alpha}^*], \operatorname{ZOH}_{\pi}[\tilde{\beta}^*])\\
            \Big\updownarrow {\scriptstyle \epsilon_4} & & \Big\updownarrow {\scriptstyle \epsilon_3}\\
            c = \tilde{J}(\tilde{\alpha}^{br}, \tilde{\beta}^{*}) & \geq & \tilde{J}(\tilde{\alpha}^*, \tilde{\beta}^*) = d.
        \end{array}
    \end{equation*}

    \underline{\textit{Step 3}.} By algebraically adding and subtracting $c$ and $d$ to the difference $b-a$, and subsequently regrouping the terms, we obtain:
    \begin{align}
        |a - b| &= b - a \nonumber \\
        &= (b - d) + (c - a) - (c - d). \label{eq:ab_two_term_exp}
    \end{align}
    Recall the fundamental property that $x \leq |x|$ for any $x \in \mathbb{R}$. Applying this inequality to the first two terms of \eqref{eq:ab_two_term_exp} yields $b - d \leq |b - d|$ and $c - a \leq |c - a| \equiv |a - c|$. Substituting these bounds into \eqref{eq:ab_two_term_exp} gives:
    \begin{equation}
        |a - b| \leq |b - d| + |a - c| - (c - d). \label{eq:ab_two_term_inter}
    \end{equation}
    Since $c \geq d$, it strictly follows that $c - d \geq 0$, rendering the term $-(c - d)$ non-positive. Removing this non-positive term from the right-hand side of \eqref{eq:ab_two_term_inter} provides a valid, looser upper bound:
    \begin{equation*}
        |a - b| \leq |a - c| + |b - d| = \epsilon_3 + \epsilon_4 + \mathcal{O}(\bar{\pi}^n).
    \end{equation*}
    This completes the proof.
\end{proof}

\section{Proof of Theorem \ref{thm:general_weak_approximation}}\label{app:thm:general_weak_approximation}
\begin{proof}
Our certification of SDG weak approximation proceeds via a two-stage procedure derived from the classical weak approximation theory of SDEs \cite{milshteinWeakApproximationSolutions1986a}: i) investigate the one-step approximation errors between the original game dynamics and the designed SDG, and ii) if the local one-step approximation errors for the first two moments are bounded by $\mathcal{O}(\bar{\pi}^2)$, then the global approximation error over the entire time horizon has an order of $1$, i.e., $\mathcal{O}(\bar{\pi})$.

Unlike classical SDE weak approximation, an SDG approximation requires both the state trajectory and the accumulated cost function to be jointly approximated, as defined in Definition \ref{def:SDG_weak_approximation}. To facilitate this, we first transform the original general game $\mathbf{G}$ and the approximated game $\mathbf{\tilde{G}}$ into auxiliary forms where the running costs are absorbed into an augmented state variable $y$.
\begin{equation*}
    \begin{aligned}
        (\mathbf{G}^{aux})\!:&\!\left\{\begin{array}{l}
             \dot{x}(t)=  f(t, x(t), u(t), v(t)), \\
             \dot{y}(t) = \mathbb{E} h(t, x(t), u(t), v(t)),\\
             x(t_0)= x_0, y(t_0) = 0,\\
             J = \mathbb{E}\left[y(T)+ g(x(T))\right].
        \end{array}\right.\\
         (\mathbf{\tilde{G}}^{aux})\!:&\!\left\{\begin{array}{l}
              d \tilde{x}(t) = \tilde{f}(t, \tilde{x}, \tilde{u}, \tilde{v}) d t + \tilde{\Sigma}^{\frac{1}{2}}(t, \tilde{x}, \tilde{u}, \tilde{v}) d W_t, \\
              d\tilde{y}(t) = \mathbb{E} \tilde{h}(t, \tilde{x}, \tilde{u}, \tilde{v}) dt,\\
              \tilde{x}(t_0) = x_0, \tilde{y}(t_0) = 0,\\
              \tilde{J} = \mathbb{E}\left[\tilde{y}(T)+ g(\tilde{x}(T))\right].
         \end{array}\right.
    \end{aligned}
\end{equation*}

We now proceed to the local one-step error analysis. Let $\delta = t_1 - t_0$ be the length of the first commitment interval, satisfying $\delta \leq \bar{\pi}$. Let $\Delta=[\Delta_x^{\top},\Delta_y]^{\top}$ be the exact one-step increment for $\mathbf{G}^{aux}$ where $\Delta_x = x(t_1)-x(t_0)$ and $\Delta_y = y(t_1)-y(t_0)$. Similarly, let $\tilde{\Delta}=[\tilde{\Delta}_x^{\top},\tilde{\Delta}_y]^{\top}$ denote the one-step increment for the approximated SDG $\mathbf{\tilde{G}}^{aux}$.

\begin{lemma}[One-step approximation]\label{lem:one_step_dg_general}
    Given that the spatial covering radii satisfy $\delta_u, \delta_v = \mathcal{O}(\bar{\pi})$, and the strategy projection maps $(\phi_1, \phi_2)$ are applied, the differences between the local moments of $\Delta$ and $\tilde{\Delta}$ conditionally on the initial state are bounded by:
    \begin{equation*}
        \begin{array}{l}
           \text{(i) } \left\| \mathbb{E} [\Delta_x] - \mathbb{E}[\tilde{\Delta}_x] \right\| = \mathcal{O}(\bar{\pi}^2), \\
           \text{(ii) } \left\| \mathbb{E} [\Delta_x \Delta_x^{\top}] - \mathbb{E}[\tilde{\Delta}_x \tilde{\Delta}_x^{\top}] \right\| = \mathcal{O}(\bar{\pi}^2), \\
            \text{(iii) } \left| \Delta_y - \tilde{\Delta}_y \right| = \mathcal{O}(\bar{\pi}^2).
        \end{array}
    \end{equation*}
\end{lemma}
\begin{proof}[Proof of Lemma \ref{lem:one_step_dg_general}]
    Under the Zero-Order Hold (ZOH) commitment pattern of $\mathbf{G}$, the continuous control processes are drawn from the probability measures $\mu$ and $\nu$ and remain constant over $[t_0, t_1)$. Expanding the exact dynamics of $\mathbf{G}^{aux}$ via Taylor expansion yields:
    \begin{equation*}
    \begin{aligned}
        \mathbb{E}[\Delta_x] &= \mathbb{E}_{\mu,\nu}[f(t_0, x_0, u, v)] \delta + \mathcal{O}(\delta^2), \\
        \mathbb{E}[\Delta_x \Delta_x^{\top}] &= \mathbb{E}_{\mu,\nu}[f(t_0, x_0, u, v) f(t_0, x_0, u, v)^{\top}] \delta^2 + \mathcal{O}(\delta^3), \\
        \Delta_y &= \mathbb{E}_{\mu,\nu}[h(t_0, x_0, u, v)] \delta + \mathcal{O}(\delta^2).
    \end{aligned}
    \end{equation*}\\
    Concurrently, the Euler-Maruyama scheme applied to the SDEs of $\mathbf{\tilde{G}}^{aux}$ yields the approximated moments. Recalling the definitions of $\tilde{f}$, $\tilde{\Sigma}$, and $\tilde{h}$ based on the discrete PMFs $(\tilde{u}, \tilde{v})$ mapped via $\phi_1, \phi_2$, we obtain:
    \begin{equation*}
    \begin{aligned}
        \mathbb{E}[\tilde{\Delta}_x] &= \tilde{f}(t_0, x_0, \tilde{u}, \tilde{v}) \delta + \mathcal{O}(\delta^2), \\
        \mathbb{E}[\tilde{\Delta}_x \tilde{\Delta}_x^{\top}] &= \tilde{\Sigma}(t_0, x_0, \tilde{u}, \tilde{v}) \delta + \tilde{f}\tilde{f}^{\top}\delta^2 + \mathcal{O}(\delta^3), \\
        \tilde{\Delta}_y &= \tilde{h}(t_0, x_0, \tilde{u}, \tilde{v}) \delta + \mathcal{O}(\delta^2).
    \end{aligned}
    \end{equation*}\\
    \underline{\textit{Proof of (i) and (iii):}} 
    Due to the Lipschitz continuity of $f$ and $h$ with respect to the control variables, the spatial quantization error introduced by the mapping to the discrete points $U_\delta$ and $V_\delta$ is bounded. Specifically, for the drift term:
    \begin{equation*}
        \left\| \mathbb{E}_{\mu,\nu}[f] - \tilde{f} \right\| \leq L_{f,u} \delta_u + L_{f,v} \delta_v = \mathcal{O}(\bar{\pi}).
    \end{equation*}
    Subtracting the expected state increments gives:
    \begin{equation*}
        \mathbb{E}[\Delta_x] - \mathbb{E}[\tilde{\Delta}_x] = \left( \mathbb{E}_{\mu,\nu}[f] - \tilde{f} \right)\delta + \mathcal{O}(\delta^2).
    \end{equation*}
    Since $\delta \leq \bar{\pi}$, the product of the spatial error $\mathcal{O}(\bar{\pi})$ and the temporal step $\delta$ yields $\mathcal{O}(\bar{\pi}^2)$. Thus, $\mathbb{E}[\Delta_x] - \mathbb{E}[\tilde{\Delta}_x] = \mathcal{O}(\bar{\pi}^2)$. The identical Lipschitz logic applies to the running cost $h$, proving (iii) that $\Delta_y - \tilde{\Delta}_y = \mathcal{O}(\bar{\pi}^2)$.\\
    \underline{\textit{Proof of (ii):}}
    For the exact dynamics of $\mathbf{G}^{aux}$, since the system state evolution $f$ is uniformly bounded over the compact action spaces, the exact second moment increment strictly scales with the squared time step:
    \begin{equation*}
        \| \mathbb{E}[\Delta_x \Delta_x^{\top}] \| = \| \mathbb{E}_{\mu,\nu}[ff^{\top}] \delta^2 + \mathcal{O}(\delta^3) \| = \mathcal{O}(\delta^2).
    \end{equation*}
    Conversely, for the approximated SDG $\mathbf{\tilde{G}}^{aux}$, the second moment increment driven by the discrete PMFs expands as:
    \begin{equation*}
        \mathbb{E}[\tilde{\Delta}_x \tilde{\Delta}_x^{\top}] = \tilde{\Sigma}(t_0, x_0, \tilde{u}, \tilde{v}) \delta + \tilde{f}\tilde{f}^{\top}\delta^2 + \mathcal{O}(\delta^3).
    \end{equation*}
    According to the preconditions established in Theorem \ref{thm:general_weak_approximation}, the designable diffusion matrix $\Sigma$, and consequently its expectation $\tilde{\Sigma}$, is positive semi-definite and strictly bounded such that $\| \tilde{\Sigma} \| \leq C_\Sigma \bar{\pi}$. Because the time step satisfies $\delta \leq \bar{\pi}$, the leading variance deviation generated by the artificial diffusion process is rigorously bounded by:
    \begin{equation*}
        \| \tilde{\Sigma} \delta \| \leq C_\Sigma \bar{\pi} \cdot \delta \leq C_\Sigma \bar{\pi}^2 = \mathcal{O}(\bar{\pi}^2).
    \end{equation*}
    Furthermore, because the expected drift $\tilde{f}$ is also bounded, the drift-induced term $\| \tilde{f}\tilde{f}^{\top} \delta^2 \|$ is bounded by $\mathcal{O}(\delta^2) \leq \mathcal{O}(\bar{\pi}^2)$. It immediately follows that the overall approximated second moment satisfies $\| \mathbb{E}[\tilde{\Delta}_x \tilde{\Delta}_x^{\top}] \| = \mathcal{O}(\bar{\pi}^2)$. \\
    By applying the triangle inequality, the difference between the exact and approximated second moments yields:
    \begin{equation*}
        \begin{aligned}
            \left\| \mathbb{E}[\Delta_x \Delta_x^{\top}] - \mathbb{E}[\tilde{\Delta}_x \tilde{\Delta}_x^{\top}] \right\| &\leq \left\| \mathbb{E}[\Delta_x \Delta_x^{\top}] \right\| + \left\| \mathbb{E}[\tilde{\Delta}_x \tilde{\Delta}_x^{\top}] \right\| \\
            &= \mathcal{O}(\delta^2) + \mathcal{O}(\bar{\pi}^2) = \mathcal{O}(\bar{\pi}^2).
        \end{aligned}
    \end{equation*}
    This confirms that the local one-step variance approximation error is theoretically preserved at order $2$, regardless of the specific bounded functional design of $\Sigma$. Hence, the proof of (ii) is completed.
\end{proof}

According to Milshtein's theorem \cite[Theorem 2]{milshteinWeakApproximationSolutions1986a}, Lemma \ref{lem:one_step_dg_general} guarantees that a local weak approximation error of $\mathcal{O}(\bar{\pi}^2)$ per step over $N = T/\bar{\pi}$ intervals accumulates to a global approximation error of $\mathcal{O}(\bar{\pi})$. Consequently, the continuous state trajectory and the accumulated cost of the designed SDG $\mathbf{\tilde{G}}^{aux}$ under pure strategies $(\phi_1(\alpha), \phi_2(\beta)) \in \tilde{\mathcal{A}}_p^\pi \times \tilde{\mathcal{B}}_p^\pi$ constitute an order $1$ weak approximation of the original game $\mathbf{G}^{aux}$ under mixed strategies $(\alpha, \beta) \in \mathcal{A}_m^\pi \times \mathcal{B}_m^\pi$. Reverting the auxiliary variables to the original cost functional yields $|J(t_0, x_0, \alpha, \beta) - \tilde{J}(t_0, x_0, \phi_1(\alpha), \phi_2(\beta))| = \mathcal{O}(\bar{\pi})$, which completes the proof.
\end{proof}

\end{appendices}

\begin{IEEEbiographynophoto}{Tao Xu}
	(S'22) received a B.S. degree in the School of Mathematical Sciences from Shanghai Jiao Tong University (SJTU), Shanghai, China. He is currently working toward the Ph.D. degree with the Department of Automation, SJTU. He is a member of Intelligent of Wireless Networking and Cooperative Control Group. His research interests include probabilistic prediction, distributionally robust optimization, dynamic games, and robotics.
\end{IEEEbiographynophoto}

\begin{IEEEbiographynophoto}{Wang Xi}
	(S'24) received the B.S. degree in the School of Electronic Information and Electrical Engineering from Shanghai Jiao Tong University (SJTU), Shanghai, China. Since 2023, he has been with the Department of Automation, SJTU. He is a member of Intelligent of Wireless Networking and Cooperative Control Group. His current research interests include differential games and its applications in mechatronic systems.
\end{IEEEbiographynophoto}

\begin{IEEEbiographynophoto}{Jianping He}
	(SM'19) is currently a full professor in the Department of Automation at Shanghai Jiao Tong University. He received the Ph.D. degree in control science and engineering from Zhejiang University, Hangzhou, China, in 2013, and had been a research fellow in the Department of Electrical and Computer Engineering at University of Victoria, Canada, from Dec. 2013 to Mar. 2017. His research interests mainly include the distributed learning, control and optimization, security and privacy in network systems.

	Dr. He serves as an Associate Editor for IEEE Trans. Control of Network Systems, IEEE Trans. on Vehicular Technology, IEEE Open Journal of Vehicular Technology, and KSII Trans. Internet and Information Systems. He was also a Guest Editor of IEEE TAC, International Journal of Robust and Nonlinear Control, etc. He was the winner of Outstanding Thesis Award, Chinese Association of Automation, 2015. He received the best paper award from IEEE WCSP'17, the best conference paper award from IEEE PESGM'17, and was a finalist for the best student paper award from IEEE ICCA'17, and the finalist best conference paper award from IEEE VTC20-FALL.
\end{IEEEbiographynophoto}

\end{document}